\documentclass[a4paper,10pt]{article}

\usepackage{amssymb}
\usepackage{amsmath}
\usepackage{enumerate}
\usepackage{graphicx}

\newtheorem{theorem}{Theorem}

\newtheorem{definition}[theorem]{Definition}
\newtheorem{proposition}[theorem]{Proposition}
\newtheorem{corollary}[theorem]{Corollary}

\newenvironment{proof}{\medskip\noindent\textbf{Proof.}}{\hspace{\stretch{1}}\rule{1ex}{1ex}\\\medskip}

\newcommand{\CQFD}{\hspace{\stretch{1}}$\Box$}
\renewcommand{\eqref}[1]{\mbox{{(\ref{#1})}}}
\newcommand{\txt}[1]{\textnormal{#1}}

\newcommand{\grandO}{\mathcal{O}}

\newcommand{\norme}[2]{\|\,#1\,\|_{#2}}
\newcommand{\hnorme}[2]{|\,#1\,|_{#2}}

\newcommand{\R}[1]{\mathbf{R}^{#1}}

\newcommand{\N}[1]{\mathbf{N}^{#1}}

\newcommand{\ii}{\txt{i}}

\newcommand{\PP}{\mathbf{P}}
\newcommand{\QQ}{\mathbf{Q}}
\newcommand{\II}{\mathbf{I}}

\newcommand{\ubeta}{u^{\beta_1,\beta_2}}
\newcommand{\ombeta}{\omega^{\beta_1,\beta_2}}
\newcommand{\FF}[1]{\omega^{F_{#1}}}
\newcommand{\uF}[1]{u^{F_{#1}}}
\newcommand{\vF}[1]{v^{F_{#1}}}
\newcommand{\rr}{\widetilde{\rho}}

\newcommand{\mpar}{m_{\shortparallel}}
\newcommand{\mopar}{m_{0,\shortparallel}}
\newcommand{\mperp}{m_{\perp}}
\newcommand{\moperp}{m_{0,\perp}}
\newcommand{\Rpar}{R_{\shortparallel}}
\newcommand{\Rperp}{R_{\perp}}
\newcommand{\SBF}{S^{LF}}
\newcommand{\SHF}{S^{HF}}
\newcommand{\Spar}{S_{\shortparallel}}
\newcommand{\SparBF}{S_{\shortparallel}^{LF}}
\newcommand{\SparHF}{S_{\shortparallel}^{HF}}
\newcommand{\SStot}{\widetilde{S}}
\newcommand{\SSBF}{\widetilde{S}^{LF}}
\newcommand{\SSHF}{\widetilde{S}^{HF}}
\newcommand{\SSpart}{\widetilde{S}_{part}}
\newcommand{\SSpar}{\widetilde{S}_{\shortparallel}}
\newcommand{\SSparBF}{\widetilde{S}_{\shortparallel}^{LF}}

\newcommand{\Xopar}{X_{0,\shortparallel}}
\newcommand{\Xpar}{X_{\shortparallel}}

\newcommand{\XNL}{X^{NL}}
\newcommand{\mupar}{\mu_{\shortparallel}}

\title{Vortex-like finite-energy asymptotic profiles for isentropic
  compressible flows}

\author{{\bf L.Miguel Rodrigues}\\
Institut Fourier, U.M.R. C.N.R.S. 5582\\ Universit\'e de Grenoble I\\
B.P. 74\\38402 Saint-Martin-d'H\`eres, France\\
{\it lmrodrig@fourier.ujf-grenoble.fr}}

\date{}

\begin{document}

\maketitle

{\hspace{\stretch{1}}} {\bf Abstract} {\hspace{\stretch{1}}}

Bidimensional incompressible viscous flows with well-localised
vorticity are well-known to develop vortex structures.  The purpose of
the present paper is to recover the asymptotic profiles describing
these phenomena for homogeneous finite-energy flows as asymptotic
profiles for near-equilibrium isentropic compressible flows. This task
is performed by extending the sharp description of the asymptotic
behaviour of near-equilibrium compressible flows obtained by David
Hoff and Kevin
Zumbrun~\cite{Hoff_Zumbrun-NS_compressible_pres_de_zero} to the case
of finite-energy \emph{vortex-like} solutions. \\

\noindent{\bf Mathematics subject classification (2000).}. 76N99,
35B40, 35M20, 35Q30.\\ 

\noindent {\bf Keywords.} Isentropic compressible Navier-Stokes
equations, near-equilibrium, long-time asymptotic profiles, vortex,
hyperbolic-parabolic composite-type, artificial viscosity
approximation.

\section*{Introduction}

The present paper is focused on the long-time asymptotic behaviour of
viscous bidimensional flows. When no exterior force is applied the
flow is expected to return to equilibrium, namely to a state of
constant density and zero velocity. Our purpose is thus to determine
asymptotic profiles for the return to equilibrium.

The motion of the considered flows may be described by the time
evolution of the pair $(\rho,u)$, $\rho=\rho\,(t,x)>0$ being the
density field of the fluid and $u=u\,(t,x)\in\R2$ the velocity
field. The main purpose of the paper is to prove that for some initial
data near equilibrium one recovers for isentropic compressible flows
the same asymptotic profiles as in the constant-density
case. Therefore let us begin introducing the constant-density profiles
we are interested in.

When the density is constant, $\rho\equiv\rho_\star$, mass
conservation and a force balance for Newtonian fluids lead to the
Navier-Stokes evolution equations
\begin{equation}\label{NS}\left.\qquad\qquad
\begin{array}{rcl}
\txt{div}\,u&=&0\\[1ex]
\partial_t\,(\rho_\star\,u)\ +\ \left(u\,\cdot\nabla\right)\,(\rho_\star\,u)
&=&\mu\,\triangle\,u\ -\ \nabla\,p
\end{array}\qquad\qquad\right\}
\end{equation}
where $\mu>0$ is the shear Lam\'e viscosity coefficient and
$p=p\,(t,x)\in\R{}$ is the pressure field of the fluid. In order to
make the former equations compatible the pressure must be determined
(up to a constant) by the elliptic equation
\begin{eqnarray}\label{incompressible}
\triangle\,p&=&-\ \rho_\star\ \txt{div}\,\big(\left(u\,\cdot\nabla\right)\,u\big)\ .
\end{eqnarray}
In this bidimensional incompressible context, it may seem more natural
and it is often more convenient to work with the curl of the velocity
rather than with the velocity itself. The evolution of the vorticity
$\omega=\partial_1\,u_2-\partial_2\,u_1$ obeys 
\begin{eqnarray}\label{NS-vorticity}
\partial_t\,\omega\ +\ u\,\cdot\,\nabla\,\omega&=&\nu\,\triangle\,\omega
\end{eqnarray}
where $\nu=\mu/\rho_\star$ is the kinematic viscosity and the velocity
$u$ is recovered by the Biot-Savart law,
\begin{eqnarray}\label{BS}
u(x)\ =\ \frac{1}{2\pi}\,\int_{\R2}\frac{(x-y)^\perp}{|x-y|^2}\ \omega\,(y)\,dy&,&x\in\R2\ ,
\end{eqnarray}
with $(z_1,z_2)^\perp=(-z_2,z_1)$, which we also denote
$u=K_{BS}\star\omega$, $K_{BS}$ being called the Biot-Savart kernel. Note
that in terms of Fourier transforms the Biot-Savart law becomes
\begin{eqnarray}\label{BS-Fourier}
\widehat{u}\,(\eta)&=&\frac{\ii\ \eta^\perp}{|\eta|^2}\ \widehat{\omega}\,(\eta)\ ,
\qquad\eta\in\R2\ .
\end{eqnarray}
Concerning the widely-developed literature about the (homogeneous)
Navier-Stokes equations, the reader is referred to some advanced
entering gates such as the following books \cite{Cannone},
\cite{Lemarie}, \cite{Lions1}, \cite{Majda_Bertozzi}, and to the more
vorticity-focused review~\cite{Ben-Artzi-review_vorticity}.

Flows with constant density and initially well-localised vorticity are
well-known to develop vortex-like structures. In the compressible case
we shall recover near equilibrium this kind of behaviour.

For instance it is proved in \cite{Gallay_W-global_stability} that any
solution $\omega$ of~\eqref{NS-vorticity} with integrable initial
datum $\omega_0$ satisfies in Lebesgue spaces
\begin{eqnarray}\label{limG}
\lim_{t\to\infty}\,t^{1-\frac1p}\ \norme{\omega(t)-\alpha\,\omega^G(t)}{p}
&=&0\ ,\\\label{limuG}
\lim_{t\to\infty}\,t^{\frac12-\frac1q}\ \norme{u(t)-\alpha\,u^G(t)}{q}
&=&0 \ ,
\end{eqnarray}
for any $1\leq p\leq\infty$ and any $2<q\leq\infty$, where
\begin{equation}\label{omG}\!
\omega^G(t,x)\ =\ \frac1t\ 
G\bigg(\frac{x}{\sqrt{\nu\,t}}\bigg)\ ,\quad
u^G(t,x)\ =\ \sqrt{\frac{\nu}{t}}\ 
v^G\bigg(\frac{x}{\sqrt{\nu\,t}}\bigg)\!
\end{equation}
with profiles
\begin{equation}\label{G}\qquad
G(\xi)\ =\ \frac{1}{4 \pi}\ e^{-|\xi|^2/4}\ ,\qquad\quad
v^G(\xi)\ =\ 
\frac{1}{2\pi}\ \frac{\xi^{\perp}}{|\xi|^2}\ \big(1-e^{-|\xi|^2/4}\big)\ ,
\end{equation}
and $\alpha\in\R{}$ is such that the initial velocity circulations at
infinity coincide at the initial time $t=0$,
\begin{eqnarray}\label{alpha}
\nu\,\alpha&=&\int_{\R2}\omega_0\,(x)\ dx\ .
\end{eqnarray}
Actually, for any $\alpha\neq0$, the vorticity $\alpha\,\omega^G$ is a
(self-similar) solution of~\eqref{NS-vorticity} with initial datum a
Dirac mass --- centred at the origin and of weight $\alpha/\nu$ ; the
corresponding flow is called Oseen vortex. Thus when the circulation
is non zero equalities~\eqref{limG} and \eqref{limuG} show that the
flow behaves asymptotically as a single vortex, whereas when
$\alpha=0$ it returns to equilibrium faster than a single vortex does.

However finite energy flows have zero circulation. Indeed, as is
easily derived from~\eqref{BS-Fourier}, to obtain an integrable
vorticity and a square-integrable velocity one must ask for the
vorticity to be of zero mean. To consider finite energy flows we must
thus investigate profiles decaying faster. Yet it is well-known,
see~\cite{Gallay_W-invariant_manifold} (combined with
\cite[Proposition 1.5]{Gallay_W-global_stability}) for instance, that
if the initial vorticity $\omega_0$ is such that
$(1+|\cdot|)^{3/2}\,\omega_0$ is square-integrable and
$\int_{\R2}\omega_0=0$ then for any index $1\leq p\leq\infty$
\begin{eqnarray}\label{limF}
\lim_{t\to\infty}\,t^{\frac32-\frac1p}\
\norme{\omega(t)-\omega^{\beta_1,\beta_2}(t)}{p}&=&0 
\end{eqnarray}
where
\begin{eqnarray}\label{ombeta}
\ombeta\,(t,x)&=&\beta_1\,\FF1\,(t,x)\ 
+\ \beta_2\,\FF2\,(t,x)
\end{eqnarray}
with for $i=1,2$
\begin{eqnarray}\label{omF}
\FF i\,(t,x)&=&\frac{1}{\sqrt{\nu}\,t^{3/2}}\
F_i\,\Big(\frac{x}{\sqrt{\nu\,t}}\Big)\\\label{F}
F_i\,(\xi)&=&\partial_i\,G\,(\xi)\ =\ -\,\frac{\xi_i}{2}\ G\,(\xi)
\end{eqnarray}
and $\beta_i$ is such that
\begin{eqnarray}\label{beta}
\nu\,\beta_i&=&-\ \int_{\R2}x_i\ \omega_0(x)\,dx\ .
\end{eqnarray}
Observe that $\ombeta$ is not a solution of
equation~\eqref{NS-vorticity} but only of its linearisation around
equilibrium, a heat equation. However equality~\eqref{limF} is easily
seen to apply also to some flows with finite measures as initial
vorticities, such as those of initial vorticity
\begin{eqnarray}
\frac{1}{2\,\nu}\ \Big(\,\delta_{(-\beta_1,0)}\ -\ \delta_{(\beta_1,0)}\,\Big)
&+&\frac{1}{2\,\nu}\ \Big(\,\delta_{(0,-\beta_2)}\ -\
\delta_{(0,\beta_2)}\,\Big)\ .
\end{eqnarray}
Moreover, defining the corresponding velocities
\begin{eqnarray}\label{ubeta}
\ubeta&=&K_{BS}\star\ombeta\ =\ \beta_1\,\uF1\,+\,\beta_2\,\uF2
\end{eqnarray}
and for $i=1,2$ 
\begin{eqnarray*}
\uF i\,(t,x)&=&(K_{BS}\star\FF i(t))\,(x)\ =\ 
\frac{1}{t}\ \vF i\,\bigg(\frac{x}{\sqrt{\nu\,t}}\bigg)\\
\vF i\,(\xi)&=&K_{BS}\star F_i\,(\xi)\ =\ \partial_i\,v^G\,(\xi)\ ,
\end{eqnarray*}
one does observe a dipole-like feature at infinity,
\begin{eqnarray*}
\vF1(\xi)&\stackrel{|\xi|\to\infty}{=}&\frac{1}{2\pi|\xi|^4}\ 
\left(\begin{array}{c}2\,\xi_1\xi_2\\\xi_2^2-\xi_1^2\end{array}\right)
+\ \grandO\,(e^{-|\xi|^2/4})\ ,\\[1em]
\vF2(\xi)&\stackrel{|\xi|\to\infty}{=}&\frac{1}{2\pi|\xi|^4}\ 
\left(\begin{array}{c}\xi_2^2-\xi_1^2\\-2\,\xi_1\xi_2\end{array}\right)
+\ \grandO\,(e^{-|\xi|^2/4})\ .
\end{eqnarray*}
Therefore equality~\eqref{limF} does show that whenever $\alpha=0$ and
$(\beta_1,\beta_2)\neq(0,0)$ the flow behaves asymptotically in time
as would do one or two pairs of vortices. Nevertheless observe
from~\eqref{BS-Fourier} that, when the vorticity $\omega$ is such that
$(1+|\cdot|)\,\omega$ is integrable (hence $\alpha$, $\beta_1$,
$\beta_2$ defined) and the velocity $u$ is integrable, parameters
$\alpha$, $\beta_1$ and $\beta_2$ must vanish and therefore the flow
should return to equilibrium again faster. From now on our attention
will be limited to these vortex-like finite-energy profiles and thus
we must eschew assuming the velocity integrable.

\begin{figure}
\begin{center}
\includegraphics[width=20em]{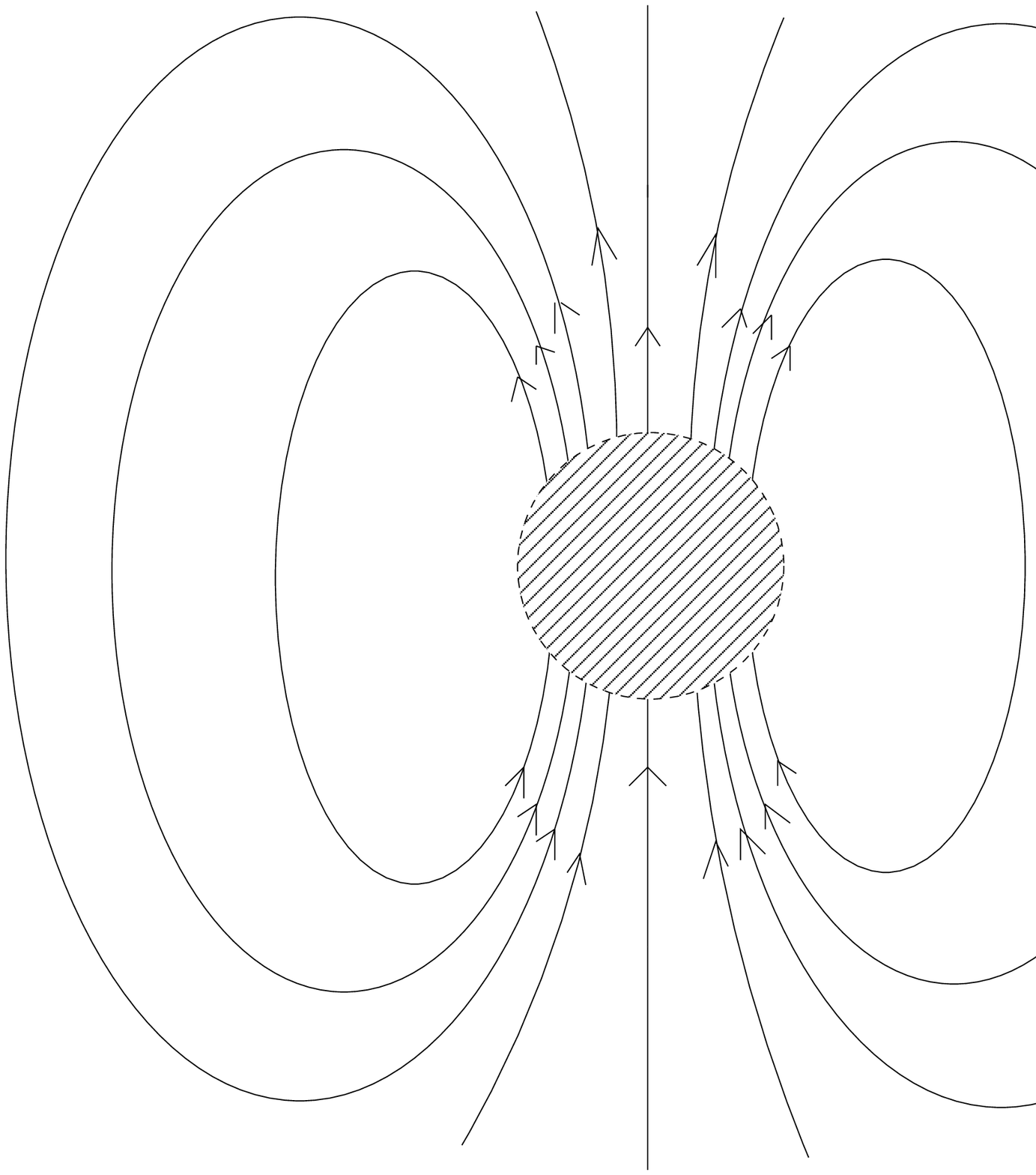}
\caption{Shape at infinity of streamlines of $\uF 1$}
\end{center}
\end{figure}


For compressible flows mass conservation and a force balance for
Newtonian fluids with constant Lam\'e coefficients give the following
equations for the time evolution of the pair $(\rho,m)$, $\rho$ being
the mass density field and $m=\rho\,u\,\in\R2$ the momentum density field,
\begin{equation}\label{NS-compressible}
\left.\begin{array}{rcl}
\partial_t\,\rho\ +\ \txt{div}\ m&=&0\\[1ex]
\partial_t\,m\ +\ \txt{div}\,(m\otimes\frac{m}{\rho})
&=&\mu\,\triangle\,(\frac{m}{\rho})
   +(\mu+\lambda)\,\nabla\,\txt{div}\,(\frac{m}{\rho})
   -\nabla\,p
\end{array}\quad\right\}
\end{equation}
where $\mu$ and $\lambda$ are the shear and bulk Lam\'e viscosity
coefficients, completed by a constitutive law for the pressure field
$p$
\begin{eqnarray}\label{law}
p&=&P\,(\rho)
\end{eqnarray}
obtained neglecting entropy variations. We require the pressure law
$P$ to be a smooth increasing function, thus the pressure increases
with density, and the Lam\'e coefficients to be such that the
viscosity tensor is elliptic namely such that $\mu>0$ and
$\lambda+2\mu>0$, which is physically relevant. Besides we choose to
formulate equations in terms of the momentum $m$ instead of the
velocity $u$ in order to keep the conservation law structure of
system~\eqref{NS-compressible}. Of course since the density $\rho$ is
expected to become asymptotically homogeneous we will also obtain
profiles for the velocity $u$.

As equations~\eqref{incompressible} and~\eqref{law} are seldom
simultaneously satisfied, there is hardly any constant-density
solution of system~\eqref{NS-compressible}. Thereby we are not
investigating stability of constant-density flows as compressible
flows, but compatibility of asymptotic behaviours for initial data
near equilibrium. To be somewhat more precise let us say that while
considering viscosity coefficients $\mu$ and $\lambda$, pressure law
$P$ and density of reference $\rho_\star$ as fixed we will ask the
density oscillations around $\rho_\star$ and the velocity to be
initially so small that both Reynolds number and Mach number shall be
small. As for compressible flows in a more general context the reader
may be referred to \cite{Majda}, \cite{Lions2} or
\cite{Feireisl-review_compressible}.

Obviously the present work is not the first one tackling the
asymptotic behaviour of near-equilibrium compressible flows. Since
1983 and the pioneer work of Kawashima about a vast class of
hyperbolic-parabolic systems of degenerate type~\cite{Kawashima-these}
equilibrium is known to be asymptotically stable for perturbations in
Sobolev spaces $H^s(\R2)$, for any integer $s$ bigger than or equal to
three. We shall make use of this stability result. Besides working
with Kawashima's solutions Hoff and Zumbrun established a precise
analysis of asymptotic behaviour of perturbations of equilibrium
\cite{Hoff_Zumbrun-NS_compressible_pres_de_zero} when initial
perturbations belong to $H^s(\R2)\cap L^1(\R2)$, for any integer $s$
bigger than five. The present paper is partially modelled on their
proof.

However the choice of Hoff and Zumbrun for initial data precludes
vortex-like asymptotic profiles as stated in~\eqref{limF}. Actually we
will obtain the same decay rates but with different profiles. These
decay rates are those of $\ombeta$ in Lebesgue spaces $L^p(\R2)$, for
$2\leq p\leq\infty$. Since these rates are not critical, decay rates
of non-linear terms should be negligible in the analysis of long-time
behaviour. An important point in the proof is that decay rates in
$L^p(\R2)$, $2\leq p\leq\infty$, are sufficient to establish that
non-linear terms can indeed be neglected in $L^q(\R2)$, for any $1\leq
q\leq\infty$, since those non-linearities are at least quadratic. Such
remark enables us to keep us away from integrability of initial data
as far as non-linear terms are concerned.

Yet in linearising around equilibrium and treating non-linear terms as
source terms we must keep in mind that
system~\eqref{NS-compressible}, \eqref{law} is quasi-linear and non
parabolic. Following~\cite{Hoff_Zumbrun-NS_compressible_pres_de_zero} we
turn round this difficulty using Kawashima's stability estimates to
bound high derivatives of the solutions. However when doing so some
terms are bounded by constants regardless of their natural decay
rates. Thus we shall require high-regularity of initial data in order
to recover natural decay rates for low derivatives of the solutions.

Let us now focus on the linearisation around $\rho=\rho_\star$ and
$m=0$ of system~\eqref{NS-compressible},
\eqref{law}. Denote the (reference) sound speed
\begin{eqnarray}
c&=&\sqrt{P'\,(\rho_\star)}
\end{eqnarray}
then consider the following system for $\rr=\rho-\rho_\star$ and $m$
\begin{equation}\label{lincompressible}
\left.\qquad\begin{array}{rcccl}
\partial_t\,\rr&+&\txt{div}\ m&=&0\\[1ex]
\partial_t\,m&+&c^2\,\nabla\rr 
&=&\mu\,\triangle\,m\ +\ (\mu+\lambda)\,\nabla\,\txt{div}\,m
\end{array}\quad.\qquad\right\}
\end{equation}
An important feature concerning system~\eqref{lincompressible} is that
it splits up into two systems, one for a curl-free part and the other
one for a constant-density divergence-free part. Let us divide
$m=\mpar+\mperp$ into its divergence free part $\mperp=\PP\,m$ and its
curl-free part $\mpar=\QQ\,m$, where $\PP$ is the Leray projection,
that is the projection onto divergence-free vector fields along
gradient fields, and $\QQ=\II-\PP$ its complementary projection. Then
system~\eqref{lincompressible} results for $(\rr,\mpar)$ in
\begin{equation}\qquad\quad\label{curlfree}\left.
\begin{array}{rcrcl}
\partial_t\,\rr&+&\txt{div}\,\mpar&=&0\\[1ex]
\partial_t\,\mpar&+&c^2\ \nabla\rr&=&(\lambda+2\mu)\ \triangle\,\mpar
\end{array}\qquad\qquad\quad\right\}
\end{equation}
and for $(0,\mperp)$ in
\begin{eqnarray}\label{divfree}
\partial_t\,\mperp\ -\ \mu\ \triangle\,\mperp&=&0\ .
\end{eqnarray}
Incidentally remark that $\mu>0$ and $\lambda+2\,\mu>0$ clearly appear to
be the conditions for ellipticity of the viscosity tensor.

Equation~\eqref{divfree} coincides with the linearisation around
equilibrium of the homogeneous Navier-Stokes equation and is thus
expected to give rise to profiles as stated in~\eqref{limF} for
suitable initial data. It remains to prove that solutions of
system~\eqref{curlfree} decay faster.

System~\eqref{curlfree} can be handled essentially as Hoff and Zumbrun
treated the whole system~\eqref{lincompressible}, the main difference
being the former system includes Leray projections in its Green
kernel.  High (and mean) frequencies of the Green kernel of
system~\eqref{curlfree} should indeed decay exponentially in time,
whereas low frequencies can be approximated by the Green kernel
$\SSpar$ of an \emph{artificial viscosity} system. This system is
derived from system~\eqref{curlfree} looking for a system whose
eigenvalues co\"incide with a second order low-frequency expansion of
eigenvalues of system~\eqref{curlfree} --- which gives a non-trivial
real part --- and that is simultaneously diagonalised with the
hyperbolic part of system~\eqref{curlfree}. This leads to a
system\footnote{See~\eqref{viscosite_artificielle} below.} of
(non-degenerate) hyperbolic-parabolic type whose hyperbolic and
parabolic parts commute. Roughly speaking, $\SSpar$'s components are
convolutions of a wave kernel and a heat kernel and look like
Gaussian functions spreading (at scale $\sqrt{(\lambda/2+\mu)\,t}\,$)
around circles scattering at scale $c\,t$ and centred at the origin.
Actually in~\cite{Hoff_Zumbrun-NS_compressible_ponctuel} the following
point-wise bounds are proved for any point $x\in\R2$ and any time $t$
bigger than one,
\begin{eqnarray}\label{ponctual}
|D^{\sigma}\SSpar(t,x)|&\leq&K\,t^{-5/4-|\sigma|/2}\ 
\left\{\begin{array}{lr}t^{3/4}\ s^{-3/2}\ ,&\qquad|x|\leq c\,(t-\sqrt t)\ ,\\
e^{-\frac{s^2}{K\,t}}\qquad\ ,&\qquad|x|\geq c\,(t-\sqrt t)\ ,\end{array}\right.
\end{eqnarray}
where $s=||x|-c\,t|$ is the distance from $x$ to the circle of radius
$c\,t$ centred at the origin. Once integrated, these bounds leads to
decay rates $t^{-(5/4-3/2p)}$ in Lebesgue space $L^p(\R2)$. Thereby
$\SSpar$ spreads faster than a heat kernel hence decays faster in
spaces requiring little localisation such as $L^p(\R2)$ for
$2<p\leq\infty$, but more slowly in $L^p(\R2)$, for $1\leq p<2$.

Let us now denote $S$ the Green kernel of
system~\eqref{lincompressible} and state the main result of this
paper, whose existence and uniqueness part is due to
Kawashima~\cite{Kawashima-these}. Lebesgue spaces are equipped with
norms $\norme{\cdot}{p}$ and Sobolev spaces $H^s(\R2)$ (based on
$L^2(\R2)$) with norms $\hnorme{\cdot}{s}$. The Green
kernel of the linearised system~\eqref{lincompressible} is denoted by
$S$.

\begin{theorem} \label{th_compressible}
Let $s$ be an integer bigger than or equal to five, $\rho_\star$ be a positive
number, $\mu$ be a positive number, $\lambda$ a real number such that
$\lambda+2\mu>0$, and $P:\R+_\star\to\R{}$ a smooth increasing function.\\
There exist positive constants $\varepsilon_0$ and $K$ and a family
$(K_p)_{1<p\leq\infty}$ of positive constants such that defining
\begin{displaymath}
X_0\ =\ (\rho_0-\rho_\star,m_0)\ ,
\qquad\Xopar\ =\ (\rho_0-\rho_\star,\mopar)\ ,
\end{displaymath}
where
\begin{displaymath}
m_0\ =\ \mopar+\moperp\ ,\qquad\mopar\ =\ \QQ\,m_0\ ,
\qquad\moperp\ =\ \PP\,m_0\ ,
\end{displaymath}
if 
\begin{displaymath}
E\ =\ \hnorme{X_0}{s}+\norme{\Xopar}{1}
+\norme{(1+|\,\cdot\,|)\ \txt{rot}\,m_0}{1}\ \leq \ \varepsilon_0\ ,
\end{displaymath}
then system~\eqref{NS-compressible}, \eqref{law} has a unique global
classical solution $(\rho,m)$ with initial datum $(\rho_0,m_0)$, and
$X=(\rho-\rho_\star,m)$ satisfies for any time $t>0$ and any
multi-index $\sigma$
\begin{enumerate}
\item when $|\sigma|\leq\frac{s-5}{2}$,
\begin{equation}\label{NL}\begin{array}{lcl}
\!\!\!\!\!\!\!\!\!\!\!\!\norme{D^{\sigma}\big(X(t)
&\!\!\!\!-\!\!\!\!&S(t)\star X_0\big)}{p}\\[1ex]
\!\!\!\!\!\!&\!\!\!\!\leq\!\!\!\!& K\,E^2\,\ln(1+t)\left\{\begin{array}{ll}
         \!\!(1+t)^{-\big(1-\frac1p+\frac{|\sigma|}{2}+\frac12\big)}
             \ \,,&2\leq p\leq\infty \\
         \!\!(1+t)^{-\big(\frac54-\frac{3}{2p}+\frac{|\sigma|}{2}+\frac12\big)}
             \,,&1\leq p\leq2\end{array}\right.\!\!;
\end{array}
\end{equation}
\item when $|\sigma|\leq\frac{s-5}{2}$, defining $\mperp=\PP\ m$,
  $\mpar=\QQ\ m$ and $\Xpar=(\rho-\rho_\star,\mpar)$,
\begin{eqnarray}
\norme{D^{\sigma}\Xpar(t)}{p}
&=&\norme{\big(D^{\sigma}\big(\rho(t)-\rho_\star\big),
          D^{\sigma}\big(m(t)-\mperp(t)\big)\big)}{p}\nonumber\\
&\leq&K\,E\ (1+t)^{-\big(\frac54-\frac{3}{2p}+\frac{|\sigma|}{2}\big)}
\ ,\qquad 2\leq p\leq\infty\ ;\label{par}
\end{eqnarray}
\item when $|\sigma|\leq\frac{s-5}{2}$, if moreover $D^\sigma\Xopar$
  is integrable then, when $t\geq1$, with
  $E'=E+\norme{D^\sigma\Xopar}{1}$,
\begin{eqnarray}\label{par2}
\norme{D^{\sigma}\Xpar(t)}{p}&\leq&
K\,E'\ t^{-\big(\frac54-\frac{3}{2p}+\frac{|\sigma|}{2}\big)}
\ ,\qquad 1\leq p<2\ ;
\end{eqnarray}
\item when $|\sigma|\leq\frac{s-5}{2}$, denoting again $\mperp=\PP\ m$,
\begin{equation}\label{perp}
\norme{D^{\sigma}\mperp(t)}{p}\ \leq\ 
K\,E\ (1+t)^{-\big(1-\frac1p+\frac{|\sigma|}{2}\big)}
\ ,\qquad 2\leq p\leq\infty\ ,
\end{equation}
and moreover, for $2\leq p\leq\infty$,
\begin{eqnarray}\label{asymperp}
\lim_{t\to\infty}\ t^{1-\frac1p+\frac{|\sigma|}{2}}\ 
\norme{D^{\sigma}\big(\mperp(t)
-\rho_\star\,\ubeta(t)\big)}{p}&=&0\ ,
\end{eqnarray}
where $\ubeta$ is defined by~\eqref{ubeta} (together
with~\eqref{ombeta}, \eqref{omF} and~\eqref{F}, remind also
$\nu=\mu/\rho_\star$) and
\begin{eqnarray}\label{betam}
\nu\,\beta_i&=&-\ \int_{\R2}x_i\
\txt{rot}\,\left(\frac{m_0}{\rho_\star}\right)(x)\,dx\ ;
\end{eqnarray}
\item if moreover $(1+|\,\cdot\,|^2)\ \txt{rot}\,m_0$ is integrable,
  then for any $1<p\leq\infty$, with $E'=E+\norme{(1+|\,\cdot\,|^2)\
  \txt{rot}\,m_0}{1}$, when $|\sigma|\leq\frac{s-5}{2}$ and $t\geq1$,
\begin{equation}\label{perp2}
\norme{D^{\sigma}\big((\,\rho(t),\ m(t)\,)
-(\,\rho_\star,\ \rho_\star\,\ubeta(t)\,)\big)}{p}
\ \leq\ K_p\,E'\,t^{-\big(\frac54-\frac{3}{2p}+\frac{|\sigma|}{2}\big)}\ .
\end{equation}
\end{enumerate}
\end{theorem}

\medskip\noindent{\bf Remarks:}\\[1ex]
\noindent{\bf1.} The hypothesis on $\txt{rot}\,m$ sufficient to
define~\eqref{betam} is enough to prove~\eqref{asymperp} yet to
quantify this asymptotic more localisation is needed, as required
for~\eqref{perp2}.\\[1ex]
\noindent{\bf2.} Estimate~\eqref{NL} does show that non-linear terms
can be neglected, whereas estimates~\eqref{par} and~\eqref{perp},
\eqref{perp2} establish that constant-density incompressible profiles
dominate in $L^p(\R2)$ for $2<p\leq\infty$ (whereas sonic waves
dominate when $1<p<2$). Indeed
\begin{eqnarray}
\lim_{t\to\infty}\ t^{1-\frac1p}\
\norme{m(t)-\rho_\star\,\ubeta(t)}{p}&=&0\ ,\qquad 2< p\leq\infty\ ,
\end{eqnarray}
while $t^{-\left(1-\frac1p\right)}$ is the decay rate of
$\rho_\star\,\ubeta(t)$ in $L^p(\R2)$ (when
$(\beta_1,\beta_2)$ is non zero) ; whereas, at least when
$(1+|\,\cdot\,|^2)\ \txt{rot}\,m_0$ is integrable,
\begin{eqnarray}
\lim_{t\to\infty}\ t^{\frac54-\frac{3}{2p}} \ \norme{(\,\rho(t),\
  m(t)\,)-\SSpar(t)\star \Xopar}{p}&=&0\ ,\qquad 1<p<2\ ,
\end{eqnarray}
where $\SSpar$ is the Green kernel of
system~\eqref{viscosite_artificielle} below, which
satisfies~\eqref{ponctual} and decays in $L^p(\R2)$ as
$t^{-\left(\frac54-\frac{3}{2p}\right)}$.
\medskip

The proof is developed in the two following sections. The next
section gathers estimates for linear equations whereas the last
one encompasses the actual proof of Theorem~\ref{th_compressible} and
in particular estimates of non-linear terms. As in
Theorem~\ref{th_compressible} from now on $\mu$, $\lambda$,
$\rho_\star$ and $P$ are considered as fixed.
\medskip

Let us also make explicit the convention used in the present paper
for Fourier transforms: when a function $f$ is integrable, its Fourier
transform is defined by
\begin{eqnarray*}
\widehat{f}\,(\eta)&=&\int_{\R2}f(x)\ e^{\ii\,\eta\cdot x}\ dx\ ,
\qquad\eta\in\R2\ .
\end{eqnarray*} At last as usual $C$ stands for a harmless constant
that may differ from line to line even in the same sequence of
inequalities. 

\section{Linear equations}

This section is devoted to the study of the system resulting from
linearisation of~\eqref{NS-compressible}, \eqref{law} around
$\rho=\rho_\star$ and $m=0$:
\begin{equation}\label{lincompressible2}
\left.\begin{array}{rcccl}
\partial_t\,\rr&+&\txt{div}\ m&=&0\\[1ex]
\partial_t\,m&+&c^2\,\nabla\rr 
&=&\mu\,\triangle\,m+(\mu+\lambda)\,\nabla\,\txt{div}\,m
\end{array}\quad\right\}
\end{equation}
where $c=\sqrt{P'(\rho_\star)}>0$ is the reference sound speed. As was
already mentioned, splitting $m=\mpar+\mperp$ into curl-free
$\mpar=\QQ\,m$ and divergence-free $\mperp=\PP\,m$ yields the system
\begin{equation}\label{linpar}\left.\qquad
\begin{array}{rcccl}
\partial_t\,\rr&+&\txt{div}\,\mpar&=&0\\[1ex]
\partial_t\,\mpar&+&c^2\,\nabla\rr&=&(\lambda+2\mu)\,\triangle\,\mpar
\end{array}\qquad\qquad\right\}
\end{equation}
and the equation
\begin{equation}\label{linperp}
\partial_t\,\mperp-\mu\,\triangle\,\mperp\ =\ 0\ .
\end{equation}

Therefore the Green kernel $S$ of system~\eqref{lincompressible2} may be
written in terms of the Green kernel $\Spar$ of system~\eqref{linpar}
and the heat kernel $K_{\mu}$, \emph{i.e.} the Green kernel of
equation~\eqref{linperp}. To make it explicit introduce kernels of
Leray projection $\PP$ and of its complementary projection $\QQ$:
\begin{displaymath}
\PP\,f\ =\ \Rperp\star f\ ,\qquad\QQ\,f\ =\ \Rpar\star f\ ,
\end{displaymath}
for any vector-field $f$. In terms of Fourier transforms note that
\begin{displaymath}
\widehat{\Rperp}(\eta)\ =\ 
\frac{\eta^{\perp}\ {}^{\txt t}\eta^{\perp}}{|\eta|^2}\ ,\qquad
\widehat{\Rpar}(\eta)\ =\ \frac{\eta\ {}^{\txt t}\eta}{|\eta|^2}\ .
\end{displaymath}
Now observe
\begin{eqnarray}\label{Sdef}
S&=&\Spar\star\left[\begin{array}{cc}\delta_0&0\\0&\Rpar\end{array}\right]
   +\left[\begin{array}{cc}0&0\\0&K_{\mu}\star\Rperp\end{array}\right]\ ,
\end{eqnarray}
where $\delta_0$ is the Dirac mass centred at the origin and of
weight one.

Keeping~\eqref{Sdef} in mind, we now study $\Spar$ and $K_{\mu}$
separately.

\subsection{Curl-free part}

Estimates of $\Spar$ used afterwards may be established as those of
$S$ developed by David Hoff and Kevin Zumbrun. That is why the needed
results shall be stated and their proofs sketched but no explicit
calculation written down, since similar calculations can be found in
\cite{Hoff_Zumbrun-NS_compressible_pres_de_zero,Hoff_Zumbrun-NS_compressible_ponctuel}.
The reader may also consult~\cite{Kobayashi_Shibata-Hoff_Zumbrun} and
references therein about some refinement for estimates of $S$ and
related subjects.

If $(\rr,\mpar)$ is a solution of system~\eqref{linpar}, the density
oscillation $\rr$ obeys
\begin{eqnarray}\label{edp_linpar}
\partial_t^2\,\rr-c^2\,\triangle\,\rr
-(\lambda+2\mu)\,\triangle\,\partial_t\,\rr&=&0\ .
\end{eqnarray}
By the way note that in the inviscid case, namely when
$\lambda=\mu=0$, the density oscillation satisfies a wave equation,
the density waves travelling at speed $c$, which is the reason why it
is called sound speed of the flow. By taking now Fourier transforms
this yields a differential equation, where $\eta$ can be thought of as
a parameter, for the quantity $y(t,\eta)=\widehat{\rr}(t,\eta)$:
\begin{eqnarray}\label{edo_linpar}
y''+(\lambda+2\mu)\,|\eta|^2\,y'+c^2\,|\eta|^2\,y&=&0\ .
\end{eqnarray}
Thereby system~\eqref{linpar} can be solved and
\begin{equation}\label{Spar_Fourier}
\widehat{\Spar}(t,\eta)\,=\,
\left[\begin{array}{cc}
\frac{\lambda^+(\eta)\,e^{\lambda^-(\eta)\,t}
-\lambda^-(\eta)\,e^{\lambda^+(\eta)\,t}}{\lambda^+(\eta)-\lambda^-(\eta)}
&-\ii\,\big(\frac{e^{\lambda^+(\eta)\,t}-e^{\lambda^-(\eta)\,t}}
{\lambda^+(\eta)-\lambda^-(\eta)}\big)\ {}^{\txt{t}}\eta\\[1em]
-\ii\,c^2\,\big(\frac{e^{\lambda^+(\eta)\,t}-e^{\lambda^-(\eta)\,t}}
{\lambda^+(\eta)-\lambda^-(\eta)}\big)\ \eta
&\frac{\lambda^+(\eta)\,e^{\lambda^+(\eta)\,t}
-\lambda^-(\eta)\,e^{\lambda^-(\eta)\,t}}{\lambda^+(\eta)-\lambda^-(\eta)}
\end{array}\right],
\end{equation}
where eigenvalues $\lambda^{\pm}$ are
\begin{eqnarray}\label{lambda+-}
\lambda^{\pm}(\eta)&=&-\frac12\,\mupar\,|\eta|^2\,
\pm\,\frac12\,\sqrt{\mupar^2\,|\eta|^4-4\,c^2\,|\eta|^2}
\end{eqnarray}
and a new viscosity parameter $\mupar$ is defined for concision's sake
by 
\begin{eqnarray}
\mupar&=&\lambda+2\,\mu\,.
\end{eqnarray}
The former formula yields also an explicit formula for $\widehat{S}$ and
enables us to perform the whole study of $\Spar$.

In order to capture the quite different behaviour of high and low
frequencies, let us split $\Spar$. Let $\chi$ be a smooth real-valued cut-off
function taking values between zero and one that is equal to one on
$\{\eta\in\R2\,|\,|\eta|\leq R_0\}$ and vanishes on
$\{\eta\in\R2\,|\,|\eta|\geq R_0+1\}$, for some $R_0>0$ to be chosen
large enough. Now divide $\Spar=\SparBF+\SparHF$ in such a way that
\begin{equation}\label{troncature}
\widehat{\SparBF}(t,\eta)\ =\ \chi(\eta)\,\widehat{\Spar}(t,\eta)\ ,
\qquad\widehat{\SparHF}(t,\eta)\ =\ (1-\chi(\eta))\,\widehat{\Spar}(t,\eta)
\end{equation}
and study separately $\SparBF$ and $\SparHF$.

\subsubsection{High frequencies}

Expending $\lambda^{\pm}(\eta)$ around $|\eta|=\infty$ gives
\begin{displaymath}
\lambda^+(\eta)\stackrel{|\eta|\to\infty}{=}
               -\frac{c^2}{\mupar}+\grandO(|\eta|^{-2})\ ,
\qquad\lambda^-(\eta)\stackrel{|\eta|\to\infty}{=}
-\mupar\,|\eta|^2+\frac{c^2}{\mupar}+\grandO(|\eta|^{-2})\ .
\end{displaymath}
\emph{A priori} high frequencies should decay exponentially. Moreover
the former expansion confirms that one component --- $\mpar$ --- should be
regularised whereas another --- $\rr$ --- should not.

To be more precise, an integral representation of
solutions of system~\eqref{linpar} using differential
equation~\eqref{edo_linpar} may be used. Define 
\begin{eqnarray*}
A(t,r)&=&\frac{1}{2\pi\ii}\int_{\mathcal{S}^+\cup\,\mathcal{S}^-}
         \frac{e^{tz}}{p\,(r,z)}\,dz\\
B(t,r)&=&\partial_t\,A(t,r)+\mupar\,r^2\,A(t,r)\\
D(t,r)&=&e^{-\mupar\,r^2t}\int_0^te^{\mupar\,r^2s}A(s,r)\,ds
\end{eqnarray*}
where $\mathcal{S}^+$ and $\mathcal{S}^-$ are circles of radius
$c^2/2\mupar$ centered at $-c^2/\mupar$ and $-\mupar\,r^2+c^2/\mupar$
respectively, and $p$ is the polynomial
$p\,(r,z)=z^2+\mupar\,r^2\,z+c^2\,r^2$. Then with obvious matricial
conventions
\begin{displaymath}
\begin{array}{rclrcl}
\widehat{\Spar}^{1,1}(t,\eta)&\!\!=\!&B(t,|\eta|)\ ,&
\widehat{\Spar}^{1,2}(t,\eta)&\!=\!&-\ii\ A(t,|\eta|)\ {}^{\txt t}\eta
\ ,\\[1ex] 
\widehat{\Spar}^{2,1}(t,\eta)&\!=&\!\!\!\!-\ii\,c^2\,A(t,|\eta|)\ \eta\ ,&
\widehat{\Spar}^{2,2}(t,\eta)
&\!\!=\!&e^{-\mupar\,|\eta|^2t}-c^2\,|\eta|^2\,D(t,|\eta|)\,.
\end{array}
\end{displaymath}
Yet expanding $1/p(z,r)$ into powers of $r^{-1}$ yields for $r$ large
enough 
\begin{eqnarray*}
A(t,r)&=&\sum_{k=0}^\infty\ A_k(t,r)\ r^{-2k-2}\\
B(t,r)&=&e^{-\frac{c^2t}{\mupar}}+\sum_{k=0}^\infty\ B_k(t,r)\ r^{-2k-2}\\
D(t,r)&=&\sum_{k=0}^\infty\ D_k(t,r)\ r^{-2k-4}
\end{eqnarray*}
with for any $k\in\N{}$
\begin{displaymath}
|A_k(t,r)|,\ |B_k(t,r)|,\ |D_k(t,r)|\ 
\leq C\ (e^{-\frac{c^2t}{\mupar}}+e^{-\frac{\mupar}{2}r^2t})\ r_0^k
\end{displaymath}
where $C$ and $r_0$ are positive constants independent of $k$, $r$ and
$t$. In quite the same way, it can also be proved that for any
$j\in\N*$
\begin{eqnarray*}
|\partial_r^jA(t,r)|
&\leq&C_j\,(e^{-\frac{c^2t}{\mupar}}+e^{-\frac{\mupar}{2}r^2t})\ r^{-j-2}\\
|\partial_r^jB(t,r)|
&\leq&C_j\,(e^{-\frac{c^2t}{\mupar}}+e^{-\frac{\mupar}{2}r^2t})\ r^{-j-2}\\
|\partial_r^jD(t,r)|
&\leq&C_j\,(e^{-\frac{c^2t}{\mupar}}+e^{-\frac{\mupar}{2}r^2t})\ r^{-j-4}\ .
\end{eqnarray*}

Now Marcinkiewic multiplier theorem (see~\cite{Stein} for instance)
and an adaptation proved
in~\cite[Proposition~4.2]{Hoff_Zumbrun-NS_compressible_pres_de_zero}
transform these estimates around $|\eta|=\infty$ into the following
proposition preceded by useful definitions.

\begin{definition}\ \\[-1em]
\begin{enumerate}[i.]
\item A bounded symbol $\widehat{f}$ is an \emph{$L^p$-multiplier} if
  the associated operator $f\star$ can be extended for any $1<p<\infty$
  from $L^2(\R2)\cap L^p(\R2)$ to $L^p(\R2)$:
\begin{displaymath}
\norme{f\star g}{p}\ \leq\ C_p\,\norme{g}{p}\ ,\qquad g\in L^p(\R2)\ ,
\qquad 1<p<\infty\ .
\end{displaymath}
\item An $L^p$-multiplier is a \emph{strong $L^p$-multiplier} if the
  above property also holds for any $1\leq p\leq\infty$:
\begin{displaymath}
\norme{f\star g}{p}\ \leq\ C\,\norme{g}{p}\ ,\qquad g\in L^p(\R2)\ ,
\qquad 1\leq p\leq\infty\ .
\end{displaymath}
\item A family of multipliers --- either strong or not --- is
  \emph{bounded} if the above constants --- $C$ or $C_p$ --- can be
  chosen uniformly for the whole family..
\end{enumerate}
\end{definition}

A typical example of $L^p$-multiplier that is not a strong
$L^p$-multiplier is given by $\widehat{\Rperp}$, the symbol of the
Leray projection $\PP$.

\begin{proposition}\label{HF}
If $R_0$ is large enough then there exists a positive
constant $b$ such that $\SparHF$, the high-frequency part of~$\Spar$,
defined by~\eqref{troncature}, satisfies
\begin{eqnarray}
\widehat{\SparHF}\,(t,\eta)&=&e^{-b\,t}\ M(t,\eta)
\end{eqnarray}
where $(M(t))_{t\geq0}$ is a bounded family of strong
$L^p$-multipliers, and for any integers $1\leq i,j,k\leq2$ with $(i,j)\neq(1,1)$
\begin{eqnarray}
\widehat{\partial_k\,\SparHF}\,{}^{i,j}\,(t,\eta)&=&
e^{-b\,t}\ (1+t^{-1/2})\ N^{i,j}_k(t,\eta)
\end{eqnarray}
where $(N^{i,j}_k(t))_{t\geq0}$ is a bounded family of
$L^p$-multipliers.
\end{proposition}

\medskip\noindent{\bf Remarks:}\\[1ex]
{\bf 1.} The only component of $\SparHF$, therefore of $\Spar$, that does
not give rise to regularisation is $(\SparHF)^{1,1}$. Indeed its
high-frequency expansion includes a Dirac mass
$e^{-c^2t/\mupar}\delta_0$. Actually since the former proposition has
been obtained bounding Fourier transforms it is still true for
$\SparHF\star\left[\begin{array}{cc}\delta_0&0\\0&\Rpar\end{array}\right]$
or even for the high-frequency part of $S$. So is it stated
in~\cite[Lemma~5.3]{Hoff_Zumbrun-NS_compressible_pres_de_zero}.\\[1ex]
{\bf 2.} Indeed components $\widehat{\partial_k\,\SparHF}\,{}^{1,2}$
and $\widehat{\partial_k\,\SparHF}\,{}^{2,1}$ are not strong
$L^p$-multipliers. Their expansions contain terms of type
$e^{-c^2t/\mupar}\ \eta_k\ \eta_{k'}/|\eta|^2$.\\[1ex]
{\bf 3.} However $\widehat{\partial_k\,\SparHF}\,{}^{2,2}$ is a strong
$L^p$-multiplier and
$\widehat{\partial_k\,\partial_{k'}\,\SparHF}\,{}^{2,2}$ is a (weak)
$L^p$-multiplier. Thereby a source term in~\eqref{linpar} would
undergo, if placed into the first equation, a regularisation of one
derivative in the determination of $\mpar$ but of none for $\rr$;
whereas, if placed in the second equation, it would undergo a
regularisation of two derivatives in the determination of $\mpar$ and
of only one for $\rr$. This is the reason why one can not establish an
existence result based on a na\"ive linearisation around equilibrium
either in momentum variables $(\rr,m)$, since $\triangle(\mpar\,\rr)$
could not be treated as a source term for the second equation, or in
velocity variables $(\rr,u)$, since $u\cdot\nabla\rr$ could not be
handled as a source term in the first equation of~\eqref{linpar}. To
turn round this difficulty
in~\cite{Danchin-NS_compressible_proche_equilibre} Rapha\"el Danchin
considers and studies a linear system including a convection term. Yet
doing so it seems difficult if not impossible to capture a precise
decay behaviour due to dispersion. We shall rather work with
Kawashima's solutions and sacrifice some regularity.

\subsubsection{Low frequencies}

From the Hausdorff-Young inequalities and explicit
formula~\eqref{Spar_Fourier} we may at once deduce the following
proposition.

\begin{proposition}\label{BF_facile}
For any multi-index $\sigma$ there exists a positive constant
$C_\sigma$ such that the low-frequency part $\SparBF$ of~$\Spar$
satisfies for any time $t\geq0$ and any real number $p$
\begin{equation}\!
\norme{D^{\sigma}\SparBF(t)}{p}\ \leq\ 
C_{\sigma}\,\left\{\begin{array}{cl}
       t^{-\big(1-\frac1p+\frac{|\sigma|}{2}\big)}&
       \txt{if}\quad\ t\geq1\quad\,\txt{and}\ 2\leq p\leq\infty\ ,\\
       1&\txt{if}\ 0\leq t\leq1\ \txt{and}\ 1\leq p\leq\infty\ .
       \end{array}\right.\!\!\!
\end{equation}
\end{proposition}

\medskip\noindent{\bf Remark:} Obviously the proposition still holds
when $S$ is substituted for $\Spar$.

\medskip
Note $\Spar$ may contain mean frequencies but since they should both
be regularised and decay exponentially the point is really in low
frequencies. Let us then perform some expansions around $\eta=(0,0)$
in order to derive a good approximation of~$\SparBF$.

As for eigenvalues we have
\begin{equation}\label{vp_en_zero}
\lambda^\pm(\eta)\ \stackrel{|\eta|\to0}{=}\ 
-\ \frac12\ \mupar\ |\eta|^2\pm\ \ii\ c\ |\eta|+\grandO(|\eta|^3)\ .
\end{equation}
We have expanded $\lambda^{\pm}$ until getting a non-trivial real part
which leads us to second order expansion. Concerning diagonalisation
basis we shall be satisfied with a first-order expansion and therefore
we look for a Green kernel diagonalised on a diagonalisation basis of
the hyperbolic part of~\eqref{linpar}. Thereby in order to build a
good low-frequencies approximation of~\eqref{linpar} we keep the same
hyperbolic part but the parabolic part is modified and we obtain
\begin{equation}\label{viscosite_artificielle}
\quad\left.\begin{array}{rcccl}
\partial_t\,\rr&+&\txt{div}\ \mpar
&=&\frac12\,\mupar\,\triangle\,\rr\\[1ex]
\partial_t\,\mpar&+&c^2\,\nabla\rr 
&=&\frac12\,\mupar\,\triangle\,\mpar\ .
\end{array}\qquad\quad\right\}
\end{equation}

As for asymptotic behaviour the Green kernel $\SSpar$ of
system~\eqref{viscosite_artificielle} should give a close
approximation of~$\Spar$. The point in the approximation is that
system~\eqref{viscosite_artificielle} is of \emph{non-degenerate}
hyperbolic-parabolic type, with hyperbolic and parabolic parts
commuting since simultaneously diagonalised.  Such a system is called
\emph{artificial viscosity system}.  See~\cite{Liu_Zeng} (where by the
way are also exposed Kawashima's estimates) to learn more about
approximations of degenerate hyperbolic-parabolic systems in the
unidimensional context
and~\cite[Section~6]{Hoff_Zumbrun-NS_compressible_pres_de_zero} for
the general case. Before establishing that~$\SSpar$ indeed
asymptotically approaches~$\Spar$, we should study the asymptotic
behaviour of $\SSpar$.

Since hyperbolic and parabolic parts of
system~\eqref{viscosite_artificielle} commute, defining $W$ the Green
kernel of hyperbolic system
\begin{equation}\label{ondes_systeme}
\left.\qquad\qquad\begin{array}{rcccl}
\partial_t\,\rr&+&\txt{div}\ \mpar&=&0\\[1ex]
\partial_t\,\mpar&+&c^2\,\nabla\rr&=&0
\end{array}\qquad\qquad\right\}
\end{equation}
 and $K_{\mupar/2}$ the heat kernel associated to
\begin{equation}
\partial_t\,f-\frac12\,\mupar\,\triangle\,f=0
\end{equation}
leads to
 $\SSpar=W\star
\left[\begin{array}{cc}K_{\mupar/2}&0\\0&K_{\mupar/2}\end{array}\right]$.
Actually since system~\eqref{ondes_systeme} implies
\begin{eqnarray}\label{ondes_equation}
\partial_t^2\ \rr\,-\,c^2\triangle\ \rr&=&0\ ,
\end{eqnarray} 
by introducing $w$ the solution of equation~\eqref{ondes_equation}
with initial datum $w(0)=0$, $\partial_tw(0)=\delta_0$, an explicit
description is obtained:
\begin{eqnarray}
\SSpar&=&\left[\begin{array}{cc}
\partial_t\,w\star K_{\mupar/2}&-\,\nabla{}^{\txt t}w\star K_{\mupar/2}\\
-c^{\,2}\,\nabla w\star K_{\mupar/2}&\partial_t\,w\star K_{\mupar/2}
\end{array}\right]\ .
\end{eqnarray}
The former formula is fully explicit since 
\begin{eqnarray}
w\,(t,x)&=&\left\{\begin{array}{lr}
\frac{1}{2\pi c}\ \frac{1}{\sqrt{c^2t^2-|x|^2}}
&\qquad\quad\txt{if}\quad|x|<c\,t\ ,\\
0&\qquad\quad\txt{if}\quad|x|\geq c\,t\ .
\end{array}\right.
\end{eqnarray}
This enables us to obtain point-wise bounds for $\SSpar$: for any
multi-index $\sigma$, there exists a positive constant $C$ such that
for any time $t\geq1$ and any point $x\in\R2$
\begin{equation}\label{ponctuel}
|D^{\sigma}\SSpar(t,x)|\ \leq\ C\,t^{-5/4-|\sigma|/2}\ 
\left\{\begin{array}{lr}
t^{3/4}\ s^{-3/2}&\ \txt{if}\quad|x|\leq c(t-\sqrt t)\ ,\\
e^{-\frac{s^2}{Ct}}&\ \txt{if}\quad|x|\geq c(t-\sqrt t)\ ,
\end{array}\right.
\end{equation}
where $s=||x|-c\,t|$ is the distance from $x$ to the circle centred at
the origin and of radius $c\,t$. The reader is referred to
\cite{Hoff_Zumbrun-NS_compressible_ponctuel} for a proof of these
estimates. Integrating in space lead then to the following
proposition. 

\begin{proposition}\label{SSparLp}
The Green kernel $\SSpar$ of the artificial viscosity
system~\eqref{viscosite_artificielle} is such that for any multi-index
$\sigma$ there exists a positive constant $C$ such that
\begin{equation}
\norme{D^\sigma\SSpar(t)}{p}\ \leq\ 
C_{\sigma}\,t^{-\big(\frac54-\frac32\frac1p+\frac{|\sigma|}{2}\big)}\ ,\qquad
t\geq1\ ,\quad 1\leq p\leq\infty\ .
\end{equation}
\end{proposition}

\medskip\noindent{\bf Remarks:}\\[1ex]
{\bf 1.} Note that combining decay rates of the heat kernel in
$L^q(\R2)$ for suitable $q$ with an estimate of the wave operator as
operator from $L^q(\R2)$ to $L^p(\R2)$ (see~\cite{Strauss}) does not
yield Proposition~\ref{SSparLp}.\\[1ex]
{\bf 2.} By getting back to~\eqref{Sdef}, since $\PP$ is not a strong
$L^p$ multiplier, it may be observed that Proposition~\ref{SSparLp} does
not give estimates of 
\begin{equation}\label{SSpart}
\SSpart\ =\ 
\SSpar\star\left[\begin{array}{cc}\delta_0&0\\0&\Rpar\end{array}\right]\ .
\end{equation}
in $L^1(\R2)$ and $L^\infty(\R2)$. However $\Rpar$ is explicit and
point-wise bounds for $\SSpart$ may indeed be obtained, leading to
\begin{equation}
\norme{D^\sigma\SSpart(t)}{p}\ \leq\ 
\frac{c_{\sigma}\,L_{\sigma}(t)}{t^{\frac54-\frac32\frac1p+\frac{|\sigma|}{2}}}
\ ,\qquad t\geq1\ ,\quad 1\leq p\leq\infty\ ,
\end{equation}
where $L_{\sigma}(t)=1+\ln t$ if $\sigma=(0,0)$ and $L_{\sigma}(t)=1$
otherwise. Once again see~\cite{Hoff_Zumbrun-NS_compressible_ponctuel}
for a proof.

We should now compare decay rates of $\Spar-\SSpar$ to the former
rates of $\SSpar$. The first part of the following proposition is
straightforward thanks to the Hausdorff-Young inequalities since $\SSpar$
is a low-frequency approximation of $\Spar$. The second part comes
from through space decomposition combining the decay rate of
$\Spar-\SSpar$ in $L^2(\R2)$ and the following point-wise bound: for
any integer $N>2$ and any multi-index $\sigma$ there exists a positive
constant $C$ such that for $t\geq1$ and $x\in\R2$, $x\neq(0,0)$
\begin{eqnarray}
|D^{\sigma}(\SSpar-\Spar)\,(t,x)|&\leq&C\,t^{-1-|\sigma|/2}\
\left(\frac{|x|}{t}\right)^{-N}\ ,
\end{eqnarray}
which is easily obtained \emph{via} a Hausdorff-Young
inequality.
See~\cite[Lemma~8.1]{Hoff_Zumbrun-NS_compressible_pres_de_zero} for a
detailed combination of these two bounds.

\begin{proposition}\label{BF}
The low-frequency part~$\SparBF$ of $\Spar$ satisfies
\begin{enumerate}
\item for any multi-index $\sigma$ there exists a positive constant
  $C_\sigma$ such that for any time $t\geq1$ and any $2\leq p\leq\infty$, 
\begin{equation}
\norme{D^{\sigma}(\SparBF(t)-\SSparBF(t))}{p}\ 
\leq\ C_{\sigma}\,t^{-\big(1-\frac1p+\frac{|\sigma|}{2}+\frac12\big)}\ ,
\end{equation}
where $\SSparBF$ is the low-frequency part of~$\SSpar$;
\item for any multi-index $\sigma$ and any real number $\theta>0$
  there exists a positive constant $C_{\sigma,\theta}$ such that for
  $t\geq1$ and $1\leq p\leq2$,
\begin{equation}
\norme{D^{\sigma}(\SparBF(t)-\SSparBF(t))}{p}\ \ \leq\ C_{\sigma,\theta}\,
t^{-\big(\frac54-\frac32\frac1p+\frac{|\sigma|}{2}+\frac12-\theta\big)}
\ ,\end{equation}
where $\SSparBF$ is the low-frequency part of $\SSpar$.
\end{enumerate}
\end{proposition}

\medskip\noindent{\bf Remark:} To prove $\SSpar$ gives a good
description of the asymptotic behaviour of $\Spar$, it only remains to
note that the high-frequency part of~$\SSpar$ also satisfies the
estimates stated in Proposition~\ref{HF} for~$\SparHF$ and it
therefore decays exponentially.

\subsection{Constant-density divergence-free part}

We now focus our attention on the linear equation for $\mperp$. Though
it is nothing but the heat equation the following estimates are not so
standard since they concern divergence-free solutions.

First note estimates for the heat kernel does not yield in a
straightforward way estimates for $K_\mu\star\Rperp$ in $L^1(\R2)$.
However once again point-wise bounds may be obtained thanks to the
Hausdorff-Young inequalities: for any multi-index $\sigma$ there
exists a positive constant $C_\sigma$ such that for any time $t>0$ and
any point $x\in\R2$
\begin{eqnarray*}
|\,D^{\sigma}(K_\mu(t)\star\Rperp)(x)\,|
&\leq&C_\sigma\ (\,\max\,(t^{1/2},|x|)\,)^{-(|\sigma|+2)}\ .
\end{eqnarray*}
See~\cite[Lemma~2.2]{Hoff_Zumbrun-NS_compressible_ponctuel} for a
proof of the former bound. Then integrating in space gives the
following proposition. (See also \cite{Fujigaki_Miyakawa} for a
different proof.)

\begin{proposition}\label{perpL1}
For any multi-index $\sigma$ that is non zero there exists a positive
constant $C_\sigma$ such that for any time $t>0$
\begin{equation}
\norme{D^{\sigma}K_\mu(t)\star\Rperp}{p}\ \leq\ 
C_{\sigma}\,t^{-\big(1-\frac1p+\frac{|\sigma|}{2}\big)}
\ ,\quad1\leq p\leq\infty\ ,\quad|\sigma|\neq0\ .
\end{equation}
\end{proposition}

The following proposition, which is the key-proposition of the present
subsection, also follows from the Hausdorff-Young inequalities. Yet since
it does not seem to be written elsewhere we shall write its proof in
full details.

\begin{proposition}\label{perpLp}\ \\[-1em]
\begin{enumerate}
\item For any multi-index $\sigma$ there exists a positive constant
  $C_{\sigma}$ such that if $\omega_0$ is such that $(1+|\,\cdot\,|)\
  \omega_0$ is integrable and $\widehat{\omega_0}(0)=0$ then the associated
  divergence-free vector-field $\moperp=K_{BS}\star\omega_0$ satisfies
  for any time $t\geq0$ 
\begin{equation}
\norme{D^{\sigma}K_{\mu}(t)\star\,\moperp}{p}\ \leq\ 
C_{\sigma}\,E\,t^{-\big(1-\frac1p+\frac{|\sigma|}{2}\big)}\ ,
\quad 2\leq p\leq\infty\ ,
\end{equation}
where $E=\norme{(1+|\,\cdot\,|)\ \omega_0}{1}$.
\item For any multi-index $\sigma$ there exists a positive constant
  $C_{\sigma}$ such that if $\omega_0$ is such that
  $(1+|\,\cdot\,|^2)\ \omega_0$ is integrable,
  $\widehat{\omega_0}(0)=0$ and
  $\nabla_{\eta}\,\widehat{\omega_0}(0)=(0,0)$ then
  $\moperp=K_{BS}\star\omega_0$ satisfies for any time $t\geq0$
\begin{equation}
\norme{D^{\sigma}K_{\mu}(t)\star\,\moperp}{p}\ \leq\ 
C_{\sigma}\,E'\,t^{-\big(1-\frac1p+\frac{|\sigma|}{2}+\frac12\big)}\ ,
\quad 2\leq p\leq\infty\ ,
\end{equation}
where $E'=\norme{(1+|\,\cdot\,|^2)\ \omega_0}{1}$
\item For any multi-index $\sigma$ there exists a positive constant
  $C_{\sigma}>0$ such that if $\omega_0$ is such that
  $(1+|\,\cdot\,|^2)\ \omega_0$ is integrable,
  $\widehat{\omega_0}(0)=0$ and
  $\nabla_{\eta}\,\widehat{\omega_0}(0)=(0,0)$ then
  $\moperp=K_{BS}\star\omega_0$ satisfies for any $t\geq0$ and any
  $2\leq p\leq\infty$ 
\begin{equation}
\norme{|\,\cdot\,|\ (D^{\sigma}K_{\mu}(t)\star\,\moperp)}{p}\ \leq\ 
C_{\sigma}\,E'\ (1+t^{-\frac12})\
t^{-\big(1-\frac1p+\frac{|\sigma|}{2}\big)}
\ ,\end{equation}
where $E'=\norme{(1+|\,\cdot\,|^2)\ \omega_0}{1}$.
\end{enumerate}
\end{proposition}

\begin{proof}
\noindent1. Since $\widehat{\omega_0}$ is Lipschitzian and
$\widehat{\omega_0}(0)=0$, $\widehat{\moperp}$ belongs to
$L^\infty(\R2)$ and for any non-zero $\eta\in\R2$
\begin{displaymath}
|\widehat{\moperp}(\eta)|\ =\ 
C\,|\eta|^{-1}\,|\widehat{\omega_0}(\eta)|\ \leq\ 
C\,\norme{\nabla_\eta\,\widehat{\omega_0}}{\infty}\ \leq\ C\,E\ .
\end{displaymath}
Yet for any $2\leq p\leq\infty$ Hausdorff-Young inequalities lead
when defining $p'$ the conjugate exponent of $p$, that is $p'$ is such
that $\frac1p+\frac{1}{p'}=1$, to
\begin{displaymath}
\norme{D^{\sigma}K_{\mu}\star\moperp}{p}\ \leq\ C\,
\norme{|\,\cdot\,|^{|\sigma|}\,e^{-\mu\,|\,\cdot\,|^2t}\,\widehat{\moperp}}{p'}
\ \leq\ C\,
t^{-\big(1-\frac1p+\frac{|\sigma|}{2}\big)}\,\norme{\widehat{\moperp}}{\infty}
\ .\end{displaymath}
Thereby the first part of the proposition is proved.\\[1ex]

\noindent2. In quite the same way from
$|\widehat{\omega_0}(\eta)|\,\leq\,C\,E'\,|\eta|^2$ we derive
\begin{displaymath}
\norme{D^{\sigma}K_{\mu}\star\moperp}{p}\ \leq\ C\ E'\ 
\norme{|\,\cdot\,|^{|\sigma|+1}\,e^{-\mu\,|\,\cdot\,|^2t}}{p'}
\ \leq\ C\,E'\,t^{-\big(1-\frac1p+\frac{|\sigma|}{2}+\frac12\big)}\ .
\end{displaymath}\\[1ex]

\noindent3. Both $|\widehat{\moperp}(\eta)|\,\leq\,C\,E'\,|\eta|$ and
\begin{displaymath}
|\nabla_{\eta}\,\widehat{\moperp}(\eta)|\ 
\leq\ C\,(\,|\eta|^{-2}\,|\widehat{\omega_0}(\eta)|
     +|\eta|^{-1}\,|\nabla_{\eta}\,\widehat{\omega_0}(\eta)|\,)\ 
\leq\ C\,E'
\end{displaymath}
stand, in such a way that
\begin{eqnarray*}
\norme{|\,\cdot\,|\ (D^{\sigma}K_{\mu}\star\moperp)}{p}&\leq&C\ 
\norme{\nabla_{\eta}\,(\widehat{D^{\sigma}K_{\mu}}\,\widehat{\moperp})}{p'}\\
&\leq&C_{\sigma}\,E'\,(\ \norme{|\,\cdot\,|^{|\sigma|}\ 
(1+t\ |\,\cdot\,|^2)\ e^{-\mu\,|\,\cdot\,|^2t}}{p'}\\
& &\qquad\qquad+\ \norme{|\,\cdot\,|^{|\sigma|}\ 
e^{-\mu\,|\,\cdot\,|^2t}}{p'}\ )\\
&\leq&C_{\sigma}\,E'\,(1+t^{-\frac12})\,
t^{-\big(1-\frac1p+\frac{|\sigma|}{2}\big)}\ ,
\end{eqnarray*}
which achieves the proof of the proposition.
\end{proof}

We now derive from the former proposition asymptotic profiles for
$\mperp$ both in $L^p(\R2)$ for $p\geq2$ without assuming further
localisation but also without obtaining convergence rates and in
$L^p(\R2)$ for $p\leq2$ when assuming more localisation for
$\txt{curl}\,\moperp$.

\begin{corollary}\label{perpLp_petit} Let $1<p\leq2$.\\
For any multi-index $\sigma$, there exists a positive constant
$C_{\sigma,p}>0$ such that if $\omega_0$ is such that
$(1+|\,\cdot\,|^2)\ \omega_0$ is integrable,
$\widehat{\omega_0}(0)=0$ and
$\nabla_{\eta}\,\widehat{\omega_0}(0)=(0,0)$ then
$\moperp=K_{BS}\star\omega_0$ satisfies for any time $t\geq0$
\begin{eqnarray}
\norme{D^{\sigma}K_{\mu}(t)\star\,\moperp}{p}\ \leq\ 
C_{\sigma,p}\ E'\ (1+t^{-\frac12})^{2\,(\frac1p-\frac12)}
\ t^{-\big(1-\frac1p+\frac{|\sigma|}{2}+\frac12\big)}\ ,\quad
\end{eqnarray}
where $E'=\norme{(1+|\,\cdot\,|^2)\ \omega_0}{1}$.
\end{corollary}

\begin{proof}
Given a non-zero function $f:\R2\to\R{}$, since $1<p\leq2$,
H\"older's inequalities yield for any $R>0$
\begin{eqnarray*}
\norme{f}{p}&\leq&\,
(\int_{|x|\leq R}|f|^p(x)\,dx\ )^{1/p}
+\,(\int_{|x|\geq R}|f|^p(x)\,dx\ )^{1/p}\\
&\leq&C_p\,\big(\,R^{\frac2p-1}\ \norme{f}{2}
+R^{\frac2p-2}\ \norme{|\,\cdot\,|\,f}{2}\,\big)
\end{eqnarray*}
hence, by choosing $R=\norme{|\,\cdot\,|\,f}{2}\ /\ \norme{f}{2}$ in order
to optimise the last term with respect to $R>0$,
\begin{eqnarray}\label{f_lemma}
\norme{f}{p}&\leq&C_p\ 
\norme{f}{2}^{2\,(1-\frac1p)}\
\norme{|\,\cdot\,|\,f}{2}^{2\,(\frac1p-\frac12)}
\end{eqnarray}
is obtained and the proof is achieved thanks to
Proposition~\ref{perpLp} applying~\eqref{f_lemma} to the function
$f=D^{\sigma}K_{\mu}\star\,\moperp$.
\end{proof}

\begin{corollary}\label{perp_asym}
  If $\omega_0$ is a real-valued function such that $(1+|\,\cdot\,|)\
  \omega_0$ is integrable, $\widehat{\omega_0}(0)=0$ and
  $\nabla_{\eta}\,\widehat{\omega_0}(0)=(0,0)$ then the associated
  divergence-free vector-field $\moperp=K_{BS}\star\omega_0$ satisfies
  for any multi-index $\sigma$
\begin{equation}
\lim_{t\to\infty}\ t^{1-\frac1p+\frac{|\sigma|}{2}}\ 
\norme{D^{\sigma}K_{\mu}(t)\star\,\moperp}{p}\ =\ 0
\ ,\qquad 2\leq p\leq\infty\ .
\end{equation}
\end{corollary}

\begin{proof}
The conclusion of the corollary has already been proved when moreover
$(1+|\,\cdot\,|^2)\ \omega_0$ is integrable. We shall obtain
Coollary~\ref{perp_asym} by a density argument. Let
$\varepsilon>0$. Choose $R_{\varepsilon}>0$ such that truncating
$\omega$ into function $\omega_{\varepsilon}$ vanishing on $\{x\,|\,|x|>
R_{\varepsilon}\}$ and coinciding with $\omega_0$ on $\{x\,|\,|x|\leq
R_{\varepsilon}\}$ yields
\begin{eqnarray*}
\norme{(1+|\,\cdot\,|)\ (\omega_0-\omega_{\varepsilon})}{1}&\leq&\varepsilon\ ,
\end{eqnarray*}
hence $|\widehat{\omega_{\varepsilon}}(0)|\leq
C\,\varepsilon$ and $|\nabla_{\eta}\,\widehat{\omega_{\varepsilon}}(0)|\leq
C\,\varepsilon$. Now define
\begin{eqnarray*}
\omega_{app}&=&\omega_{\varepsilon}
-[\,\widehat{\omega_{\varepsilon}}(0)]\ G
-\ii\,[\partial_{\eta_1}\,\widehat{\omega_{\varepsilon}}(0)]\ F_1
-\ii\,[\partial_{\eta_2}\,\widehat{\omega_{\varepsilon}}(0)]\ F_2
\end{eqnarray*}
to obtain a function $\omega_{app}$ localised as a Gaussian function satisfying
$\widehat{\omega_{app}}(0)=0$,
$\nabla_{\eta}\,\widehat{\omega_{app}}(0)=(0,0)$ and
\begin{eqnarray*}
\norme{(1+|\,\cdot\,|)\ (\omega_0-\omega_{app})}{1}&\leq&C\,\varepsilon\ .
\end{eqnarray*}
Let $2\leq p\leq\infty$ and $\sigma$ be a multi-index. The second part
of Proposition~\ref{perpLp} yields a $t_{\varepsilon}>0$ such that
for $t\geq t_{\varepsilon}$
\begin{eqnarray*}
  t^{1-\frac1p+\frac{|\sigma|}{2}}\ 
  \norme{D^{\sigma}K_{\mu}(t)\star\,K_{BS}\star\,\omega_{app}}{p}
  &\leq&\varepsilon\ .\end{eqnarray*}
Then, with a constant independent of $\varepsilon$, the triangle inequality
together with the first part of Proposition~\ref{perpLp} give for any
$t\geq t_{\varepsilon}$
\begin{eqnarray*}
t^{1-\frac1p+\frac{|\sigma|}{2}}\ 
\norme{D^{\sigma}K_{\mu}(t)\star\,\moperp}{p}&\leq&C\,\varepsilon\ ,
\end{eqnarray*}
which achieves the proof and the present section.
\end{proof}

\section{Non-linear terms}

Now taking advantage of estimates for the Green kernel $S$ of the
linearised system~\eqref{lincompressible} we prove
Theorem~\ref{th_compressible}. It only remains to bound non-linear
terms. This task shall be performed in two steps. First we establish
estimates of $X(t)=(\rho(t)-\rho_\star,m(t))$ and non-linear terms in
Lebesgue spaces $L^p(\R2)$, for $2\leq p\leq\infty$, by a continuity
fix-point-like argument. Then we use these bounds to estimates
non-linear terms in $L^p(\R2)$, for $1\leq p<2$.

For the sake of conciseness, write
\begin{eqnarray*}
X(t)&=&S(t)\star\,X_0\ +\ \XNL(t)\ .
\end{eqnarray*}
Then
\begin{eqnarray}\label{nonlin_compressible}
\XNL(t)&=&\sum_{k=1}^2\ \int_0^tS(t-t')\star\partial_k\,Q_k(t')\,dt'
\end{eqnarray}
where, for $k=1,2$,
\begin{displaymath}
Q_k\ =\ Q_k^1+Q_k^2\ ,\qquad 
Q_k^1\ =\ \left(\begin{array}{c}0\\q_k^1\end{array}\right)\ ,\qquad
Q_k^2\ =\ \sum_{k'=1}^2\left(\begin{array}{c}0\\
                 \partial_{k'}\,q_k^{2,k'}\end{array}\right)\ ,
\end{displaymath}
in such a way that
\begin{eqnarray*}
\!\!\!\sum_{k=1}^2\ \partial_k\,q_k^1
&=&-\,\txt{div}\,\Big(\ m\otimes\frac{m}{1+\rr}\ \Big)
-\nabla\,\big(\,P\,(1+\rr)-c^2\,\rr\,\big)\\
\!\!\!\sum_{k,k'=1}^2\partial_k\,\partial_{k'}\,q_k^{2,k'}\!\!\!
&=&-\ \mu\ \triangle\,\Big(\ \frac{m\ \rr}{1+\rr}\ \Big)\ 
-\ (\mu+\lambda)\ \nabla\,\txt{div}\,\Big(\ \frac{m\ \rr}{1+\rr}\ \Big)\ .
\end{eqnarray*}

\subsection{The case $p\geq2$}

As was already mentioned, the Green kernel $S$ does not regularise
enough to enable us to deal ingenuously with all terms arising in
$\XNL$ and we resort to Kawashima's estimates
\cite{Kawashima-these}. Again note that doing so we bound some
quantities by constants regardless of their natural decay rates.

\begin{theorem}[Kawashima, 1983 \cite{Kawashima-these}]\label{Kawashima}
Let $s\geq3$ be an integer.\\
There exist $\varepsilon_0>0$ and $C>0$ such that if $X_0=(\rr_0,m_0)$
belongs to $H^s(\R2)$ with
\begin{displaymath}
E\ =\ \hnorme{X_0}{s}\ \leq\ \varepsilon_0
\end{displaymath}
then system~\eqref{nonlin_compressible} has a unique global classical
solution $X=(\rr,m)$ of initial datum $X_0$, satisfying for any time
$t\geq0$
\begin{eqnarray*}
\hnorme{X(t)}{s}^2\ +\ \int_0^t\hnorme{\nabla\,X(t')}{s-1}^2\ dt'
&\leq&C\,E^2\ .
\end{eqnarray*}
\end{theorem}

As in Theorem~\ref{th_compressible} we assume $s\geq5$ and prove, when 
\begin{eqnarray*}
E&=&\hnorme{X_0}{s}\ 
+\ \norme{\Xopar}{1}\ +\ \norme{(1+|\,\cdot\,|)\,\txt{rot}\ m_0}{1}
\end{eqnarray*}
is small, that for any $2\leq p\leq\infty$ and any multi-index $\sigma$
such that $|\sigma|\leq s-4$
\begin{eqnarray}\label{NL_p>2}
\norme{D^{\sigma}X(t)}{p}&\!\!\leq\!\!&CE\ 
(1+t)^{-\big(1-\frac1p+\frac12\min(|\sigma|,\,s-4-|\sigma|)\big)}\\\!
\norme{D^{\sigma}\XNL(t)}{p}&\!\!\leq\!\!&
CE^2\,\ln(1\!+\!t)\ 
(1\!+\!t)^{-\big(1-\frac1p+\frac12\min(|\sigma|,\,s-5-|\sigma|)+\frac12\big)}
\ .\end{eqnarray}

For this purpose,
following~\cite{Hoff_Zumbrun-NS_compressible_pres_de_zero}, we introduce
\begin{eqnarray*}
A(t)&=&\sup_{\substack{0\leq t'\leq t\\[0.5ex]
2\leq p\leq\infty\\[0.5ex]|\sigma|\leq s-4}}
(1+t')^{1-\frac1p+\frac12\min(|\sigma|,\,s-4-|\sigma|)}
\ \norme{D^{\sigma}X(t')}{p}\\
B(t)&=&\sup_{\substack{0<t'\leq t\\[0.5ex]
2\leq p\leq\infty\\[0.5ex]|\sigma|\leq s-4}} 
\frac{(1+t')^{1-\frac1p+\frac12\min(|\sigma|,\,s-5-|\sigma|)+\frac12}}
{\ln(1+t')}
\,\norme{D^{\sigma}\XNL(t)}{p}\ .\end{eqnarray*}
The present subsection is essentially devoted to the proof of the
following inequality
\begin{eqnarray}\label{B(t)}
B(t)&\leq&C\ (\,E^2+A(t)^2+A(t)^{s-2}\,)\ .
\end{eqnarray}
Together with linear estimates it shall yield
\begin{eqnarray*}
A(t)&\leq&C\ (\,E+A(t)^2+A(t)^{s-2}\,)
\end{eqnarray*}
enabling us to propagate, whenever $2\,C\,E<1$ and $4\,C^2E<1/2$,
both $A(t)+A(t)^{s-3}\leq1/2$ and
\begin{displaymath}
A(t)\ \leq\ \frac{C\,E}{1-A(t)-A(t)^{s-3}}\ \leq\ 2\,C\,E\ ,
\end{displaymath}
which may be plugged in~\eqref{B(t)}. Therefore as for the purpose of
the present subsection it is enough to prove~\eqref{B(t)}.

In order to establish~\eqref{B(t)} divide $S$ into a low-frequency
part $\SBF$ and a high-frequency part $\SHF$ as was done for $\Spar$
in~\eqref{troncature} and split $\XNL$ into
\begin{eqnarray}
\XNL(t)&=&\sum_{k=1}^2\ \int_0^{t/2}\SBF(t-t')\star\partial_k\,Q_k(t')\,dt'
\nonumber\\
&+&\sum_{k=1}^2\ \int_{t/2}^t\SBF(t-t')\star\partial_k\,Q_k^1(t')\,dt'
\nonumber\\
&+&\!\!\sum_{k,k'=1}^2
\int_{t/2}^t\SBF(t-t')\star\partial_k\,\partial_{k'}\,Q_k^{2,k'}(t')\,dt'
\nonumber\\
&+&\sum_{k=1}^2\ \int_0^t\SHF(t-t')\star\partial_k\,Q_k(t')\,dt'
\nonumber\\
&=&\XNL_1(t)+\XNL_2(t)+\XNL_3(t)+\XNL_4(t)\ .
\end{eqnarray}
Let $2\leq p\leq\infty$ and $\sigma$ a multi-index such that
$|\sigma|\leq s-4$.\\[1ex]

\medskip\noindent{\bf 1.} For some multi-indices $\sigma'$ of length
$|\sigma'|=|\sigma|+1$ Young's inequality yields
\begin{eqnarray*}
\norme{D^{\sigma}\XNL_1(t)}{p}&\leq&
\sum_{\sigma'}\int_0^{t/2}\norme{D^{\sigma'}\SBF(t-t')\star\,Q(t')}{p}\,dt'\\
&\leq&C\,\sum_{\sigma'}
\int_0^{t/2}\norme{D^{\sigma'}\SBF(t-t')}{p}\ \norme{Q(t')}{1}\,dt'\\
&\leq&C\,(1+t/2)^{-\big(1-\frac1p+\frac{|\sigma|}{2}+\frac12\big)}\,
       \int_0^{t/2}\norme{Q(t')}{1}\,dt'\ .
\end{eqnarray*}
As $X(t)$ is bounded in $L^{\infty}(\R2)$ thanks to Sobolev'
embeddings and Theorem~\ref{Kawashima}, from Theorem~\ref{Kawashima}
may be derived
\begin{eqnarray}\label{Q(t)L1}
\int_0^t\norme{Q(t')}{1}\,dt'&\leq&C\,
\int_0^t(\,\norme{X(t')}{2}^2+\norme{\nabla\,X(t')}{2}^2\,)\ dt'
\nonumber\\
&\leq&C\,
\int_0^t(\,A(t')^2(1+t')^{-1}+\norme{\nabla\,X(t')}{2}^2\,)\ dt'
\nonumber\\
&\leq&C\,(\,E^2+A(t)^2)\ \ln(1+t)
\end{eqnarray}
by using, when $0\leq t\leq1$, $\int_0^t\norme{\nabla\,X}{2}^2\leq
C\,E^2\,t$. Thereby
\begin{equation}
\norme{D^{\sigma}\XNL_1(t)}{p}\ \leq\ C\,(E^2+A(t)^2)\ \ln(1+t)
\ (1+t)^{-\big(1-\frac1p+\frac{|\sigma|}{2}+\frac12\big)}.
\end{equation}

\medskip\noindent{\bf 2.} When defining $1\leq r\leq2$ by
$1+1/p=1/2+1/r$, Young's inequality yields
\begin{eqnarray*}
\norme{D^{\sigma}\XNL_2(t)}{p}&\leq&\sum_{k=1}^2
\int_{t/2}^t\norme{\partial_k\,\SBF(t-t')\star D^{\sigma}Q_k^1(t')}{p}\,dt'\\
&\leq&C\,\int_{t/2}^t\norme{\nabla\,\SBF(t-t')}{2}\  
                    \norme{D^{\sigma}Q^1(t')}{r}\,dt'\\
&\leq&C\,\int_{t/2}^t(1+t-t')^{-1}\ \norme{D^{\sigma}Q^1(t')}{r}\,dt'\ .
\end{eqnarray*}
Now for such an $r$, H\"older's inequalities combined with Leibniz'
rule for differentiation give
\begin{eqnarray*}\!
\norme{D^{\sigma}Q^1(t)}{r}
&\!\leq\!&C\,\sum_{\substack{\sum|\sigma_i|=|\sigma|}}
\norme{D^{\sigma_1}X(t)}{p}\ \norme{D^{\sigma_2}X(t)}{2}\ 
\prod_{i\geq3}\,\norme{D^{\sigma_i}X(t)}{\infty}\\
&\!\leq\!&C\,(A(t)^2+A(t)^{\max(|\sigma|,2)})\\
&\!\!&\times\ (1+t)^{-\big(1-\frac1p+\frac12+
\underset{\sum|\sigma_i|=|\alpha|}{\min}
\underset{i}{\sum}\frac12\min(|\sigma_i|,\,s-4-|\sigma_i|)\big)}\\
&\!\leq\!&C\,(A(t)^2+A(t)^{\max(|\sigma|,2)})\ 
(1\!+\!t)^{-\big(1-\frac1p+\frac12\min(|\sigma|,\,s-4-|\sigma|)+\frac12\big)}
\!,\end{eqnarray*}
the last inequality being proved in the following way: consider
$\sigma_i$'s such that $\sum|\sigma_i|=|\alpha|$, then either all
$\sigma_i$'s satisfies $|\sigma_i|\leq (s-4)/2$ and 
\begin{displaymath}
\sum_i\ \min(|\sigma_i|,\,s-4-|\sigma_i|)\ 
=\ \sum_i\ |\sigma_i|\ =\ |\sigma|\ ,
\end{displaymath}
or there is a $\sigma_{i_0}$ such that $|\sigma_{i_0}|>(s-4)/2$ hence
\begin{displaymath}
\sum_i\min(|\sigma_i|,\,s-4-|\sigma_i|)\geq 
\min(|\sigma_{i_0}|,\,s-4-|\sigma_{i_0}|)=s-4-|\sigma_{i_0}|
\geq s-4-|\sigma|\,.
\end{displaymath}
Therefore
\begin{eqnarray}
\norme{D^{\sigma}\XNL_2(t)}{p}&\leq&C\,
(A(t)^2+A(t)^{\max(|\sigma|,2)})\,\ln(1+t)
\nonumber\\
&&\times\ (1+t)^{-\big(1-\frac1p+\frac12\min(|\sigma|,\,s-4-|\sigma|)
+\frac12\big)}.
\end{eqnarray}

\medskip\noindent{\bf 3.} Assume first that $\sigma$ is
non-zero. Again letting $1\leq r\leq2$ be such that $1+1/p=1/2+1/r$
Young's inequality yields for some $\sigma'$ such that
$|\sigma'|=|\sigma|-1$
\begin{eqnarray*}
\norme{D^{\sigma}\XNL_3(t)}{p}&\leq&C\ \sum_{\sigma'}\ 
\int_{t/2}^t\norme{D^2\,\SBF(t-t')}{2}\ \norme{D^{\sigma'}Q^2(t')}{r}\,dt'\\
&\leq&C\ \sum_{\sigma'}\ 
\int_{t/2}^t(1+t-t')^{-3/2}\ \norme{D^{\sigma'}Q^2(t')}{r}\,dt'\ .
\end{eqnarray*}
Now again
\begin{displaymath}
\norme{D^{\sigma'}Q^2(t')}{r}\ \leq\ C\,(A(t)^2+A(t)^{\max(|\sigma|,2)})\,
(1+t)^{-\big(1-\frac1p+\frac12\min(|\sigma|,\,s-4-|\sigma|)+\frac12\big)}.
\end{displaymath}
Therefore when $|\sigma|\neq0$
\begin{eqnarray}
\norme{D^{\sigma}\XNL_3(t)}{p}&\leq&C\,
(A(t)^2+A(t)^{\max(|\sigma|,2)})\,\min\,(1,t)\nonumber\\
&&\times\ (1+t)^{-\big(1-\frac1p+\frac12\min(|\sigma|,\,s-4-|\sigma|)
+\frac12\big)}.
\end{eqnarray}
Putting one derivative less on $\SBF$ it may also be proved that
\begin{eqnarray}
\norme{\XNL_3(t)}{p}&\leq&C\,A(t)^2\,\ln(1+t)
\ (1+t)^{-\big(1-\frac1p+\frac12\min(1,\ s-5)+\frac12\big)}.\qquad
\end{eqnarray}

\medskip\noindent{\bf 4.} For some $\sigma'$ such that $|\sigma'|=|\sigma|+1$
\begin{eqnarray*}
\norme{D^{\sigma}\XNL_4(t)}{p}&\leq&C\ \sum_{\sigma'}\ 
\int_0^te^{-b(t-t')}\,\norme{D^{\sigma'}Q(t')}{p}\ dt'\ .
\end{eqnarray*}
Now on the one hand when $|\sigma'|\leq (s-4)$
\begin{eqnarray*}
\norme{D^{\sigma'}Q^1(t)}{p}&\leq&C\,(A(t)^2+A(t)^{\max(|\sigma|+1,2)})\\
&&\times\ (1+t)^{-\big(1-\frac1p+\frac12\min(|\sigma|+1,\,s-5-|\sigma|)+1\big)}
\end{eqnarray*}
and when $|\sigma'|\leq (s-5)$
\begin{eqnarray*}
\norme{D^{\sigma'}Q^2(t)}{p}&\leq&C\,(A(t)^2+A(t)^{|\sigma|+2})\\
&&\times\
(1+t)^{-\big(1-\frac1p+\frac12\min(|\sigma|+2,\,s-6-|\sigma|)+1\big)}.
\end{eqnarray*}
And on the other hand when $|\sigma'|=(s-3)$ Sobolev embeddings and
Theorem~\ref{Kawashima} yield
\begin{eqnarray*}
\norme{D^{\sigma'}Q(t)}{p}&\leq&C\ A(t)\ (E+A(t)+A(t)^{|\sigma|+1})\ 
(1+t)^{-\big(1-\frac1p\big)}.
\end{eqnarray*}
At last, when $|\sigma'|=(l-2)$, $\norme{D^{\sigma'}Q^2(t)}{p}$ may be
bounded by interpolation in such a way that in any case stands
\begin{eqnarray}\label{DQ(t)}
\norme{D^{\sigma'}Q(t)}{p}&\leq&C\,(\,E^2+A(t)^2+A(t)^{|\sigma|+2})\nonumber\\
&&\times\ (1+t)^{-\big(1-\frac1p+\frac12\min(|\sigma|,\,s-5-|\sigma|)
+\frac12\big)}.
\end{eqnarray}
It should be emphasised that the former estimates are the critical
ones leading to some loss of decay in the estimates. Observe now that
for $\gamma>1$
\begin{eqnarray*}
\int_0^{t/2}e^{-b(t-t')}\,(1+t')^{-\gamma}\,dt'&\leq&C\,e^{-bt/2}\\
\int_{t/2}^te^{-b(t-t')}\,(1+t')^{-\gamma}\,dt'&\leq&C\,(1+t)^{-\gamma}.
\end{eqnarray*}
Therefore gathering everything leads to
\begin{eqnarray}
\norme{D^{\sigma}\XNL_4(t)}{p}&\leq&C\,(\,E^2+A(t)^2+A(t)^{|\sigma|+2})
\nonumber\\
&&\times\ (1+t)^{-\big(1-\frac1p+\frac12\min(|\sigma|,\,s-5-|\sigma|)
+\frac12\big)}.
\end{eqnarray}

This achieves the proof of~\eqref{B(t)}.\CQFD

\subsection{The case $p<2$}

It only remains to bound $\XNL(t)$ and its derivatives in $L^p(\R2)$, for
$1\leq p\leq2$. An important point is that we shall only make use of bounds
of $X(t)$ in $L^p(\R2)$ for $2\leq p\leq\infty$ thus we do not need to
assume $X(t)$ to be integrable !

Our aim is to prove for any multi-index $\sigma$ such that
$|\sigma|\leq(l-2)$ and any index $1\leq p\leq2$
\begin{equation}\label{NL_p<2}
\norme{D^\sigma\XNL(t)}{p}\,\leq\,C\,E^2\,\ln(1+t)\ 
(1+t)^{-\big(\frac54-\frac32\frac1p+\frac12\min(|\sigma|,\,s-5-|\sigma|)
+\frac12\big)}\ .\!
\end{equation} Since $(\widehat{S}(t))_{0\leq t\leq2}$ is a family of
bounded strong $L^p$-multipliers and Kawashima's theorem gives through
Sobolev' embeddings $\norme{D^{\sigma}Q(t)}{p}\leq\,C\,E^2$ whenever
$|\sigma|\leq(s-4)$, then does stand
\begin{eqnarray*}
\norme{D^\sigma\XNL(t)}{p}&\leq&C\,E^2\,t\ ,\qquad0\leq t\leq2\,,
\quad1\leq p\leq\infty\ .
\end{eqnarray*}
The point is therefore in dealing with $\XNL(t)$ and its derivatives
for $t\geq2$.

For this purpose, having in mind~\eqref{Sdef}, we introduce $\SStot$
defined by
\begin{eqnarray}\label{SSdef}
\SStot&=&
\SSpar\star\left[\begin{array}{cc}\delta_0&0\\0&\Rpar\end{array}\right]
+\left[\begin{array}{cc}0&0\\0&K_{\mu}\star\Rperp\end{array}\right]\ ,
\end{eqnarray}
where $\SSpar$ is the Green kernel of the artificial viscosity
system~\eqref{viscosite_artificielle}. Then $\SStot$ is the Green
kernel of the artificial viscosity system
\begin{equation}\label{SS_systeme}
\left.\begin{array}{rcccl}
\partial_t\,\rr&+&\txt{div}\ m&=&(\mu+\frac{\lambda}{2})\ \triangle\,\rr\\[1ex]
\partial_t\,m&+&c^2\,\nabla\rr 
&=&\mu\ \triangle\,m\ +\ \frac{\lambda}{2}\ \nabla\ \txt{div}\,m
\end{array}\qquad\right\}
\end{equation}
and should well approach $S$. Indeed Proposition~\ref{BF} still holds
when replacing $\Spar$ with $S$ and $\SSpar$ with $\SStot$. Now split $S$
and $\SStot$ into high-frequency and low-frequency parts and divide
$\XNL$ into
\begin{eqnarray}
\XNL(t)&=&\sum_{k=1}^2\ \int_{t-1}^tS(t-t')\star\partial_k\,Q_k(t')\,dt'
\nonumber\\
&+&\sum_{k=1}^2\ \int_0^{t/2}\SStot(t-t')\star\partial_k\,Q_k(t')\,dt'
\nonumber\\
&+&\sum_{k=1}^2\ \int_{t/2}^{t-1}\SStot(t-t')\star\partial_k\,Q_k^1(t')\,dt'
\nonumber\\
&+&\!\!\sum_{k,k'=1}^2
\int_{t/2}^{t-1}\SStot(t-t')\star\partial_k\,\partial_{k'}\,Q_k^{2,k'}(t')\,dt'
\nonumber\\
&+&\sum_{k=1}^2\ \int_0^{t/2}(\SBF-\SSBF)(t-t')\star\partial_k\,Q_k(t')\,dt'
\nonumber\\
&+&\sum_{k=1}^2\ 
\int_{t/2}^{t-1}(\SBF-\SSBF)(t-t')\star\partial_k\,Q_k^1(t')\,dt'
\nonumber\\
&+&\!\!\sum_{k,k'=1}^2\int_{t/2}^{t-1}(\SBF-\SSBF)(t-t')
\star\partial_k\,\partial_{k'}\,Q_k^{2,k'}(t')\,dt'
\nonumber\\
&+&\sum_{k=1}^2\ \int_0^{t-1}(\SHF-\SSHF)(t-t')\star\partial_k\,Q_k(t')\,dt'
\nonumber\\
&=&\XNL_1(t)+\cdots+\XNL_8(t)\ .
\end{eqnarray}
Let $1\leq p\leq2$ and $\sigma$ be such that
$|\sigma|\leq(s-4)$.\\[1ex]

\medskip\noindent{\bf 1.} Since $(\widehat{S}(t))_{0\leq t\leq1}$ is a
bounded family of strong $L^p$-multipliers, we may obtain
\begin{eqnarray}
\norme{D^{\sigma}\XNL_1(t)}{p}
&\leq&C\,\int_{t-1}^t\norme{\nabla D^{\sigma}Q(t')}{p}\,dt'\nonumber\\
&\leq&C\ E^2\,\int_{t-1}^t
(1+t')^{-\big(1-\frac1p+\min(|\sigma|,\,s-5-|\sigma|)+\frac12\big)}\,dt'
\nonumber\\
&\leq&C\ E^2\ t^{-\big(1-\frac1p+\min(|\sigma|,\,s-5-|\sigma|)+\frac12\big)}\ ,
\end{eqnarray}
where $\norme{D^{\sigma}Q(t')}{p}$ is bounded mainly as was
established~\eqref{DQ(t)} in the former subsection.\\[1ex]

\noindent{\bf 2.} When $t\geq2$, for some $\sigma'$ such that
$|\sigma'|=|\sigma|+1$, a Young inequality yields
\begin{eqnarray}
\norme{D^{\sigma}\XNL_2(t)}{p}&\leq&C\ \sum_{\sigma'}\ 
\int_0^{t/2}\norme{D^{\sigma'}\SStot(t-t')}{p}\ \norme{Q(t')}{1}\,dt'
\nonumber\\ 
&\leq&C\ t^{-\big(\frac54-\frac32\frac1p+\frac{|\sigma|}{2}+\frac12\big)}\,
\int_0^{t/2}\norme{Q(t')}{1}\,dt'\nonumber\\
&\leq&C\,E^2\ \ln(1+t)\
t^{-\big(\frac54-\frac32\frac1p+\frac{|\sigma|}{2}+\frac12\big)}\ ,
\end{eqnarray}
when using~\eqref{Q(t)L1} with $A(t)\leq C\,E$.\\[1ex]

\noindent{\bf 3.} Since $5/4-3/2p+1/2\leq1$, H\"older's and Young's
inequalities yield through a change of variables
\begin{eqnarray}
\norme{D^{\sigma}\XNL_3(t)}{p}&\leq&C\,
\int_{t/2}^{t-1}\norme{\nabla\SStot(t-t')}{p}
\ \norme{D^{\sigma}Q^1(t')}{1}\,dt'\nonumber\\ 
&\leq&CE^2\,\int_{t/2}^{t-1}
\Big[(t-t')^{-\big(\frac54-\frac32\frac1p+\frac12\big)}
\nonumber\\[-1em]
&&\qquad\qquad\qquad\qquad\quad
\times(1+t')^{-\big(1+\frac12\min(|\sigma|,\,s-4-|\sigma|)\big)}\Big]\,dt'
\nonumber\\
&\leq&CE^2\ t^{-\big(\frac54-\frac32\frac1p+\frac12
+\frac12\min(|\sigma|,\,s-4-|\sigma|)\big)}\ 
\int_{1-\frac1t}^{1/2}\frac{dt''}{t''^{\frac54-\frac32\frac1p+\frac12}}
\nonumber\\
&\leq&CE^2\,\ln(1+t)\ \ 
t^{-\big(\frac54-\frac32\frac1p+\frac12\min(|\sigma|,\,s-4-|\sigma|)
+\frac12\big)},
\end{eqnarray}
where $\norme{D^{\sigma}Q^1(t')}{1}$ has been estimated mainly as was
$\norme{D^{\sigma}Q^1(t')}{p}$ in the former subsection.\\[1ex]

\noindent{\bf 4.} When $\sigma$ is non-zero acting in quite the same
way leads for some multi-indices $\sigma'$ such that
$|\sigma'|=|\sigma|-1$ to
\begin{eqnarray}
\norme{D^{\sigma}\XNL_4(t)}{p}&\leq&C\,\int_{t/2}^{t-1}
\norme{D^2\SStot(t-t')}{p}\ \norme{D^{\sigma'}Q^2(t')}{1}\,dt'\nonumber\\ 
&\leq&CE^2\,\int_{t/2}^{t-1}\Big[(t-t')^{-\big(\frac54-\frac32\frac1p+1\big)}
\nonumber\\[-1em]
&&\qquad\qquad\qquad\qquad\quad
\times(1+t')^{-\big(1+\frac12\min(|\sigma|,\,s-4-|\sigma|)\big)}\Big]\,dt'
\nonumber\\
&\leq&CE^2\,\int_{t/2}^{t-1}\Big[
(t-t')^{-\big(\frac54-\frac32\frac1p+\frac12\big)}
\nonumber\\[-1em]
&&\qquad\qquad\qquad\qquad\quad
\times(1+t')^{-\big(1+\frac12\min(|\sigma|,\,s-4-|\sigma|)\big)}\Big]\,dt'
\nonumber\\
&\leq&CE^2\,\ln(1+t)\ \ 
t^{-\big(\frac54-\frac32\frac1p+\frac12\min(|\sigma|,\,s-4-|\sigma|)
+\frac12\big)},
\end{eqnarray}
where has been used $(t-t')^{-\frac12}\leq1$ in the
integrand. In a similar way, since $5/4-3/2p+1/2\leq1/2$,
\begin{eqnarray}
\norme{\XNL_4(t)}{p}&\leq&
C\,\int_{t/2}^{t-1}\norme{\nabla\SStot(t-t')}{p}\ \norme{Q^2(t')}{1}\,dt'
\nonumber\\ 
&\leq&CE^2\,\ln(1+t)\ \ 
t^{-\big(\frac54-\frac32\frac1p+\frac12\min(1,\,s-5)+\frac12\big)}.
\end{eqnarray}

\noindent{\bf 5.} By proceeding as for $\XNL_2$ may be obtained when
$t\geq2$
\begin{eqnarray}
\norme{D^{\sigma}\XNL_5(t)}{p}
&\leq&C\,E^2\,\ln(1+t)\ \ 
t^{-\big(\frac54-\frac32\frac1p+\frac{|\sigma|}{2}+1-\theta\big)},
\end{eqnarray}
for some $0<\theta\leq1/2$.\\[1ex]

\noindent{\bf 6.} Proceeding as for $\XNL_3$ and taking into account
$(t-t')^{-(1/2-\theta)}\leq1$, when $0<\theta\leq1/2$, in the
integrand lead to
\begin{equation}
\norme{D^{\sigma}\XNL_6(t)}{p}\ \leq\ C\ E^2\ \ln(1+t)\ \ 
t^{-\big(\frac54-\frac32\frac1p+\frac12\min(|\sigma|,\,s-4-|\sigma|)
+\frac12\big)}.
\end{equation}

\noindent{\bf 7.} Proceeding as for $\XNL_4$ and taking into account
$(t-t')^{-(1-\theta)}\leq1$, whenever $0<\theta\leq1/2$, in the
integrand give
\begin{equation}
\norme{D^{\sigma}\XNL_7(t)}{p}\ \leq\ C\,E^2\,\ln(1+t)\ \ 
t^{-\big(\frac54-\frac32\frac1p+\frac12\min(|\sigma|,\,l-2-|\sigma|)
+\frac12\big)}.
\end{equation}

\noindent{\bf 8.} The high-frequency study yields
\begin{eqnarray*}
\norme{D^{\sigma}\XNL_8(t)}{p}&\leq
&C\,\int_0^{t-1}e^{-b(t-t')}\,\norme{\nabla D^{\sigma}Q(t')}{p}\,dt'\\
&\leq&C\,E^2\int_0^te^{-b(t-t')}\,
(1+t')^{-\big(1-\frac1p+\min(|\sigma|,\,s-5-|\sigma|)+\frac12\big)}\,dt'.
\end{eqnarray*}
Now observe when $\gamma\geq0$ for any time $t\geq0$ 
\begin{eqnarray*}
\int_0^{t/2}e^{-b(t-t')}\,(1+t')^{-\gamma}\,dt'&\leq&C\ t\ e^{-bt/2}\ ,\\
\int_{t/2}^te^{-b(t-t')}\,(1+t')^{-\gamma}\,dt'&\leq&C\ (1+t)^{-\gamma}.
\end{eqnarray*}
Thereby
\begin{eqnarray}
\norme{D^{\sigma}\XNL_8(t)}{p}&\leq&C\ E^2\ \ 
(1+t)^{-\big(1-\frac1p+\min(|\sigma|,\,s-5-|\sigma|)+\frac12\big)}.\qquad
\end{eqnarray}

This achieves the proof of Theorem~\ref{th_compressible}.
Estimate~\eqref{NL} comes from~\eqref{NL_p>2} and~\eqref{NL_p<2}.
Estimates~\eqref{par} and~\eqref{par2} are consequences of
estimate~\eqref{NL} and linear estimates of Propositions \ref{HF},
\ref{BF_facile}, \ref{SSparLp} \& \ref{BF}. Estimate~\eqref{par} may
be derived from estimate~\eqref{NL} and estimates of Propositions
\ref{perpL1} \& \ref{perpLp}. Equality~\eqref{asymperp} is deduced from
estimate~\eqref{NL} and Corollary~\ref{perp_asym}. At last
estimate~\eqref{perp2} may be obtained from estimate~\eqref{NL} and
Proposition~\ref{perpLp_petit}.\CQFD

\bigskip\bigskip

\noindent {\bf Acknowledgements.} Again I warmly thank Thierry Gallay
for having supported me along this work.

\bibliographystyle{plain} 
\bibliography{Ref}

\end{document}